\documentclass[12pt,a4paper]{amsart}

\usepackage{amssymb}  
\usepackage{latexsym} 

\newcommand{\game}[2]{\left\{#1 \mid #2\right\}}
\newcommand{\gameL}[1]{\left\{#1\,\,\right|}
\newcommand{\gameR}[1]{\left.\,\,#1\right\}}
\newcommand{\gameBig}[2]{\Bigl\{#1 \Bigm| #2\Bigr\}}

\newcommand{\fu}{\mathrel{\,\rule[-1pt]{0.3pt}{9pt}\,\rule[-1pt]{0.3pt}{9pt}\,}}
\newcommand{\lf}{\lhd}  
\newcommand{\gf}{\rhd}  

\newcommand{\up}{\mathord{\uparrow}}
\newcommand{\down}{\mathord{\downarrow}}
\renewcommand{\star}{\mathord{\ast}}

\newcommand{\mex}{\operatorname{mex}}
\renewcommand{\And}{\mathrel{\text{\quad and\quad}}}

\newcommand{\Z}{{\mathbb Z}}
\newcommand{\N}{{\mathbb N}}
\newcommand{\Q}{{\mathbb Q}}
\newcommand{\R}{{\mathbb R}}
\newcommand{\C}{{\mathbb C}}

\newtheorem{Theorem}{Theorem}[section]
\newtheorem{Proposition}[Theorem]{Proposition}
\newtheorem{Lemma}[Theorem]{Lemma}
\newtheorem{Corollary}[Theorem]{Corollary}

\theoremstyle{definition}

\newtheorem{Definition}[Theorem]{Definition}
\newtheorem{Remark}[Theorem]{Remark}

\numberwithin{equation}{section}

\newenvironment{Proof}{\par\noindent{\sc Proof:}}%
                         {\hspace*{\fill}\nobreak$\Box$\par\medskip}

\newcommand{\Strut}{\raisebox{-6pt}{\rule{0pt}{17pt}}}
\newcommand{\lineclear}{\hfill\\}

\newcommand\hide[1]{}

\def\footnotenonum#1{
\def\thefootnote{\fnsymbol{footnote}}
\footnotetext[0]{\hskip -4 ex #1}
\def\thefootnote{\arabic{footnote}}
}
\def\MSC#1{\footnotenonum{{\bf MSC Classification (2000): }#1}}
\def\keywords#1{\footnotenonum{{\bf Keywords: }#1}}

\begin{document}

\author{Dierk Schleicher}

\address{School of Engineering and Science,
            International University Bremen,
            Postfach 750~561, D-28725 Bremen, Germany.}
\email{dierk@iu-bremen.de}

\author{Michael Stoll}

\address{School of Engineering and Science,
            International University Bremen,
            Postfach 750~561, D-28725 Bremen, Germany.
}
\email{m.stoll@iu-bremen.de}

\title{An Introduction to Conway's Games and Numbers}

\date{\today}
\MSC{91-02, 91A05, 91A46, 91A70}

\maketitle


\section{Combinatorial Game Theory}
\label{Sect:CGT}

Combinatorial Game Theory is a fascinating and rich theory, based on a
simple and intuitive recursive definition of games, which yields
a very rich algebraic structure: games can be added and
subtracted in a very natural way, forming an abelian GROUP
(\S~\ref{Sect:Group}). There is a distinguished sub-GROUP of games called
{\em numbers} which can also be multiplied and which form a FIELD
(\S~\ref{Sect:Field}): this field contains both the real numbers
(\S~\ref{Ssect:Short}) and the ordinal numbers (\S~\ref{Sect:Ordinal}) (in
fact, Conway's definition generalizes both Dedekind sections and von
Neumann's definition of ordinal numbers). All Conway numbers can be
interpreted as games which can actually be played in a natural way; in
some sense, if a game is identified as a number, then it is understood
well enough so that it would be boring to actually play it
(\S~\ref{Sect:GamesNumbers}). Conway's theory is deeply satisfying from a
theoretical point of view, and at the same time it has useful
applications to specific games such as Go \cite{Go}. There is a beautiful
microcosmos of numbers and games which are infinitesimally close to zero
(\S~\ref{Sect:Infini}), and the theory contains the classical and
complete Sprague-Grundy theory on impartial games
(\S~\ref{Sect:Impartial}).

The theory was founded by John H.\ Conway in the 1970's. Classical
references are the wonderful books {\em On Numbers and Games} \cite{ONAG}
by Conway, and {\em Winning Ways} by Berlekamp, Conway and Guy \cite{WW};
they have recently appeared in their second editions.
\cite{WW} is a most beautiful book bursting with examples and results but
with less stress on mathematical rigor and exactness of some statements.
\cite{ONAG} is still the definitive source of the theory, but rather
difficult to read for novices; even the second edition shows that it was
originally written in one week, and we feel that the order of presentation
(first  numbers, then games) makes it harder to read and adds unnecessary
complexity to the exposition. \cite{SN} is an entertaining story about
discovering surreal numbers on an island.

This note attempts to furnish an introduction to Combinatorial Game
Theory that is easily accessible and yet mathematically precise and
self-contained, and which provides complete statements and proofs for
some of the folklore in the subject.
We have written this note with readers in mind who have enjoyed looking at
books like \cite{WW} and are now eager to come to terms with the
underlying mathematics, before embarking on a deeper study in
\cite{ONAG,GONC} or elsewhere. While this note should be complete enough
for readers without previous experience with combinatorial game theory,
we recommend looking at \cite{WW}, \cite{GONC} or \cite{AGBB} to pick up
the playful spirit of the theory. We felt no need for duplicating many
motivating examples from these sources, and we have no claims for
originality on any of the results.

{\sc Acknowledgements}. This note grew out of a summer school on the subject
which we taught for students of the {\em Studienstiftung des deutschen Volkes}
in 2001 (in La Villa). We would like to thank the participants for the
inspiring and interesting discussions and questions which eventually led to
the writing of this note. We would also like to thank our Russian friends, in
particular Alexei Belov, for encouragement and interest in Conway games.


\section{The GROUP of Games}
\label{Sect:Group}


\subsection{What is a game?} \label{SSect:Game}
\cite[\S\S~7, 0]{ONAG}, \cite[\S\S~1, 2]{WW}

Our notion of a game tries to formalize the abstract structure underlying
games such as Chess or Go: these are two-person games with complete
information; there is no chance or shuffling. The two players are usually
called {\em Left} ($L$) and {\em Right} ($R$). Every game has
some number of {\em positions}, each of which is described by the set of
allowed {\em moves} for each player. Each move (of Left or Right) leads
to a new position, which is called a (left or right) {\em option} of the
previous position. Each of these options can be thought of as another
game in its own right: it is described by the sets of allowed moves for
both players.

   From a mathematical point of view, all that matters are the sets of left
and right options that can be reached from any given position --- we can
imagine the game represented by a rooted tree with vertices representing
positions and with oriented edges labeled $L$ or $R$ according to the
player whose moves they reflect. The root represents the initial
position, and the edges from any position lead to another rooted
(sub-)tree, the root of which represents the position just reached.

Identifying a game with its initial position, it is completely described
by the sets of left and right options, each of which is another game.
This leads to the recursive Definition~\ref{DefGame} (\ref{DefGame1}).
Note that the sets $L$ and $R$ of options may well be infinite or empty.
The Descending Game Condition (\ref{DefGame2}) simply says that every game
must eventually come to an end no matter how it is played; the number of
moves until the end can usually not be bounded uniformly in terms of the
game only.

\begin{Definition}[Game] \label{DefGame}
     \begin{enumerate}
       \item \label{DefGame1}
         Let $L$ and $R$ be two sets of games. Then the ordered pair
         $G:=(L,R)$ is a {\em game.}
       \item \label{DefGame2}
         (Descending Game Condition (DGC)). There is no infinite sequence of
games
         $G^i=(L^i,R^i)$ with $G^{i+1}\in L^i\cup R^i$ for all $i\in\N$.
     \end{enumerate}
\end{Definition}

Logically speaking, this recursive definition does not tell you what games
{\em are}, and it does not need to: it only needs to specify the
axiomatic properties of games. A major purpose of this paper is of course to
explain the meaning of the theory; see for example the creation of games
below.

\begin{Definition}[Options and Positions] \strut
     \begin{enumerate}
       \item
         (Options). The elements of $L$ and $R$ are called {\em left} resp.\
         {\em right options} of $G$.
       \item
         (Positions). The {\em positions} of~$G$ are~$G$ and all the positions
         of any option of~$G$.
     \end{enumerate}
\end{Definition}

In the recursive definition of games, a game consists of two sets of
games. Before any game is `created', the only set of games we have is the
empty set: the simplest game is the `zero game' $0 = (\{\;\},\{\;\})$
with $L=R=\{\;\}$: in this game, no player has a move.
Now that we have a non-empty set of games, the next simpler games are
$1 = \left(\{0\},\{\;\}\right)$ (whose name indicates that it represents
one free move for Left), $-1 = \left(\{\;\}, \{0\}\right)$ (a free move for
Right) and $\star = \left(\{0\},\{0\}\right)$ (a free move for whoever
gets to take it first).

{\bf Notation}.
We simplify (or abuse?) notation as follows: let
$L=\{G^{L_1},G^{L_2},\dots\}$ and $R=\{G^{R_1},G^{R_2},\dots\}$ be two
arbitrary sets of games (we do {\em not} mean to indicate that $L$ or $R$
are countable or non-empty); then for
\[
G=(L,R)=\left(\{G^{L_1},G^{L_2},\dots\},\{G^{R_1},G^{R_2},\dots\}\right)
\]
we write $G=\game{G^{L_1},G^{L_2},\dots}{G^{R_1},G^{R_2},\dots}$. Hence a
game is really a set with two distinguished kinds of elements: the left
respectively right options%
\footnote{
It is customary to abuse notation and write $\game{L}{R}$ for the ordered
pair $(L,R)$. We will try to avoid that in this paper.}.
With this notation, the four simplest games
introduced so far can be written more easily as
\[
0=\game{\;}{\;} \qquad
1=\game{0}{\;} \quad
-1=\game{\;}{0} \qquad
\star=\game{0}{0} \,\,.
\]

Eventually, we will want the two players to move
alternately: that will be formalized in \S~\ref{Ssec:Winning}; but the
Descending Game Condition will be needed to hold even when players do not
move alternately, see \S~\ref{Ssec:Adding}.

The simple recursive (and at first mind-boggling) definition of games
has its counterpart in the following equally simple induction principle
that is used in almost every proof in the theory.

\begin{Theorem}[Conway Induction] \label{ThmConwayInduction}
     Let $P$ be a property which games might have, such that any game $G$ has
     property $P$ whenever all left and right options of $G$ have this
     property. Then every game has property $P$.

     More generally, for $n\geq 1$, let $P(G_1,\dots,G_n)$ be a
     property which any $n$-tuple of games might have (i.e., an $n$-place
     relation). Suppose that $P(G_1,\dots,G_i,\dots,G_n)$ holds whenever all
     $P(G_1,\dots,G'_i,\dots,G_n)$ hold (for all $i\in\{1,\dots,n\}$ and all
     left and right options $G'_i\in L_i\cup R_i$, where $G_i=(L_i,R_i)$). Then
     $P(G_1,\dots,G_n)$ holds for every $n$-tuple of games.
\end{Theorem}
\begin{Proof}
     Suppose there is a game $G$ which does not satisfy $P$. If all left and
     right options of $G$ satisfy $P$, then $G$ does also by hypothesis,
     so there is an option $G'$ of $G$ which does not satisfy $P$.
     Continuing this argument inductively, we obtain a sequence
     $G,G',G'',\dots$ of games, each an option of its predecessor,
     which violates the Descending Game Condition.
     Note that formalizing this argument needs the
     Axiom of Choice.

     The general statement follows similarly: if
     $P(G_1,\dots,G_i,\dots,G_n)$ is false, it follows that some
     $P(G_1,\dots,G'_i,\dots,G_n)$ (for some $i$ and some $G'_i\in L_i\cup
     R_i$) is also false, so either some $P(G_1,\dots,G''_i,\dots,G_n)$ or some
     $P(G_1,\dots,G'_i,\dots,G'_j,\dots,G_n)$ is false, and it is easy to
     extract an infinite sequence of games $G_i,G'_i,G''_i,\dots$ which are
     options of their predecessors, again a contradiction.
\end{Proof}

Note that Conway Induction does not need an explicit induction base (as
opposed to ordinary induction for natural numbers which must be based at
$0$): the empty game $0=\game{\;}{\;}$ satisfies property $P$ automatically
because all its options do---there is no option which might fail
property $P$.

As a typical illustration of how Conway Induction works, we show that its
first form implies the Descending Game Condition.

\begin{Proposition}
     Conway Induction implies the Descending Game Condition.
\end{Proposition}
\begin{Proof}
     Consider the property $P(G)$: {\em there is no infinite chain of games
     $G,G',G'',\dots$ starting with $G$ so that every game is followed by one
     of its options}. This property clearly is of the kind described by
     Conway Induction, so it holds for every game.
\end{Proof}

Conway's definition of a game~\cite[\S\S~0, 7]{ONAG} consists of
part~(\ref{DefGame1})
in Definition~\ref{DefGame}, together with the statement {\em `all games
are constructed in this way'}. One way of making this precise is by Conway
Induction: a game is `constructed in this way' if all its options are,
so Conway's axiom becomes a property which all games enjoy. We have chosen
to use the equivalent Descending Game Condition in the definition in order
to treat induction for one or several games on an equal footing.

Another easy consequence of the Conway Induction principle is that
the positions of a game form a set (and not a proper CLASS).


\subsection{Winning a game} \label{Ssec:Winning}
\cite[\S\S~7, 0]{ONAG}, \cite[\S~2]{WW}

  From now on, suppose that the two players must move alternately.
When we play a game, the most important aspect usually is whether we
can win it or will lose it. In fact, most of the theory is about deciding
which player can force a win in certain kinds of games. So we need some
formal definition of who wins and who loses; there are no ties or draws in
this theory\footnote{This is one reason why Chess does not fit well into
our theory; another one is that addition is not natural for Chess. The
game of Go, however, fits quite well.}. The basic decision we make here is
that we consider a player to have lost a game when it is his turn to move
but he is unable to do so (because his set of options is empty): the idea
is that we cannot win if we do not have a good move, let alone no
move at all. This {\em Normal Play Convention}, as we will see, leads to a
very rich and appealing theory.

There is also a {\em Mis\`ere Play Convention} that the loser is the one
who makes the last move; with that convention, most of our theory would
fail, and there is no comparably rich theory known: our fundamental
equality $G=G$ for every game $G$ very much rests on the normal play
convention (Theorem~\ref{ThmFundamentalInequality}); see also the end of
Section~\ref{Sect:Impartial} for the special case of impartial games.
Another possible winning convention is by score; while scores are not
built into our theory, they can often be simulated: see the remark after
Definition~\ref{DefLeftRightStops} and \cite[Part~3]{WW}.

Every game $G$ will be of one of the following four outcome classes:
(1) Left can enforce a win, no matter who starts; (2) Right can enforce a
win, no matter who starts; (3) the first player can enforce a win, no
matter who it is; (4) the second player can enforce a win, no matter who.
We will abbreviate these four possibilities by $G > 0$, $G < 0$, $G \fu
0$, and $G = 0$, respectively: here, $G\fu 0$ is usually read `$G$ is
fuzzy to zero'; the justification for the notation $G=0$ will become clear
in \S~\ref{Ssec:Adding}. We can contract these as usual: $G\ge 0$ means
$G>0$ or $G=0$, i.e.\ Left can enforce a win (at least) if he is the
second player; $G\le 0$ means that Right can win as second player;
similarly, $G\gf 0$ means $G>0$ or $G\fu 0$, i.e.\ Left can win as first
player (`$G$ is greater than or fuzzy to zero'), and $G\lf 0$ means that
Right can win as first player.

It turns out that only $G\ge 0$ and $G\le 0$ are fundamental: if $G\ge
0$, then Left wins as second player, so Right has no good opening move. A
good opening move for Right would be an option $G^R$ in which Right could
win; since Left must start in $G^R$, this would mean $G^R\le 0$. This
leads to the following formal definition:

\begin{Definition}[Order of Games] \label{DefOrder1}
     We define:
     \begin{itemize}
       \item
         $G\ge 0$ unless there is a right option $G^R\le 0$;
       \item
         $G\le 0$ unless there is a left option $G^L\ge 0$;
     \end{itemize}
\end{Definition}
The interpretation of winning needs to be based at games where Left or
Right win immediately: this is the {\em Normal Play Convention} that a
player loses when it is her turn but she has no move available. Formally,
if Left has no move at all in $G$, then clearly $G\le 0$ by definition, so
Right wins when Left must start but cannot move.
Note that the convention `both players move alternately' enters the
formal theory in Definition~\ref{DefOrder1}.

As so often, Definition~\ref{DefOrder1} is recursive: in order to decide
whether or not $G\ge 0$, we must know whether $G^R\le 0$ etc. The
Descending Game Condition makes this well-defined: if there was a game
$G$ for which $G\ge 0$ or $G\le 0$ was not well-defined, then this could
only be so because there was an option $G^L$ or $G^R$ for which these
relations were not well-defined etc., and this would eventually violate
the DGC.

It is convenient to introduce the following conventions.
\begin{Definition}[Order of Games] \label{DefOrder2}
     We define:
     \begin{itemize}
       \item
         $G=0$ if $G\ge 0$ and $G\le 0$, i.e.\ there are no options 
$G^R\le 0$ or
         $G^L\ge 0$;
       \item
         $G>0$ if $G\ge 0$ but not $G\le 0$, i.e.\ there is an option $G^L\ge 0$
         but no $G^R\le 0$;
       \item
         $G<0$ if $G\le 0$ but not $G\ge 0$, i.e.\ there is an option $G^R\le 0$
         but no $G^L\ge 0$;
       \item
         $G\fu 0$ if neither $G\ge 0$ nor $G\le 0$, i.e.\ there are options
         $G^L\ge 0$ and $G^R\le 0$.
       \item
         $G\gf 0$ if $G\le 0$ is false, i.e.\ there is a left option $G^L\ge 0$;
       \item
         $G\lf 0$ if $G\ge 0$ is false, i.e.\ there is a right option 
$G^R\le 0$;
     \end{itemize}
\end{Definition}

A game $G$ such that $G = 0$ is often called a `zero game' (not to be
confused with {\em the} zero game $0 = \game{\;}{\;}$!).

All these cases can be interpreted in terms
of winning games; for example, $G\gf 0$ means that Left can win when
moving first: indeed, the condition assures the existence of a good
opening move for Left to $G^L\ge 0$ in which Left plays second.

Note that these definitions immediately imply the claim made above that
for every game $G$ exactly one of the following statements
is true: $G = 0$, $G > 0$, $G < 0$ or $G \fu 0$. They are the four cases
depending on the two independent possibilities $\exists G^L\ge 0$ and
$\exists G^R\le 0$, see also Figure~\ref{TabOutcome}.

\begin{figure}[htbp]
     \begin{center}
       \begin{tabular}{|rl|c|c|}
         \hline
                               &           &
         \multicolumn{2}{c|}{if Right starts, then} \\
                               &           & Left wins & Right wins \\
         \hline
         if Left starts, & Left wins & $G>0$ & $G\fu 0$\Strut \\ \cline{2-4}
         then                                                   & Right wins
         & $G=0$ & $G<0$\Strut \\
         \hline
       \end{tabular}
     \end{center}
     \caption{The four outcome classes.}
     \label{TabOutcome}
\end{figure}

When we say `Left wins' etc., we mean that Left
can enforce a win by optimal play; this does {\em not} mean that we
assume a winning strategy to be actually known, or that Right might not
win if Left plays badly. For example, for the beautiful game of
Hex~\cite{Hex1,Hex2}, there is a simple proof that the first player can
enforce a win, though no winning strategy is known unless the board size
is very small --- and there are serious Hex tournaments.%
\footnote{The rules are usually modified to eliminate the first player's
advantage. With the modified rules, one can prove that the second player
can enforce a win (if he only knew how!), and the situation is then similar.}

The existence of a strategy
for exactly one player (supposing that it is fixed who starts) is built
into the definitions: to fix ideas, suppose that Right starts in a game
$G\ge 0$. Then there is no right option $G^R\le 0$, so either $G$ has no
right option at all (and Left wins effortlessly), or all $G^R\gf 0$, so
every $G^R$ has a left option $G^{RL}\ge 0$: whatever Right's
move, Left has an answer leading to another game $G^{RL}\ge 0$, so Left
will never be the one who runs out of moves when it is his turn. By the
Descending Game Condition, the game eventually stops at a position where
there are no options left for the player whose turn it is, and this must
be Right. Therefore, Left wins. Note that this argument does not assume
that a strategy for Left is known, nor does it provide an explicit
strategy.

We will define equality of games below as an equivalence relation. What we
have so far is equality of games in the set-theoretic sense; to distinguish
notation, we use a different word for this and say that two games $G$ and
$H$ are {\em identical} and write $G \equiv H$ if they have the same sets
of  (identical) left resp.\ right options.

For the four simplest games, we have the following outcome classes.
We have obviously $0 = 0$, since no player has a move; then
it is easy to see that $1 > 0$, $-1 < 0$ and $\star \fu 0$.

{\bf A note on set theory.}
Definition~\ref{DefGame} might look simple and innocent, but the CLASS of
games thus defined is a proper CLASS (as opposed to a set): one way to
see this is to observe that every ordinal number is a game
(\S~\ref{Ssect:Ordinal}). We have adopted the convention (introduced
in~\cite{ONAG}) of writing GROUP, FIELD etc.\ for algebraic structures
that are proper classes  (as opposed to sets).

The set-theoretic foundations of our theory are the Zermelo-Fraenkel
axioms, including the axiom of choice (ZFC), and expanded by proper
classes. It is a little cumbersome to  express our theory in terms of
ZFC: a game is a set with two kinds of elements, and it might be more
convenient to treat combinatorial game theory as an appropriately
modified analog to ZFC. See the discussion in \cite[Appendix to Part
0]{ONAG}, where Conway argues that ``the complicated nature of these
constructions [expressing our theory in terms of ZFC] tells us more about
the nature of formalizations within ZF than about our system of numbers
\dots [formalization within ZFC] destroys a lot of its symmetry''. In this
note, we will not go into details concerning such issues; we only note
that our Descending Game Condition in Definition~\ref{DefGame}
corresponds to the Axiom of Foundation in ZFC.
In the special case of {\em impartial games}, however,
the Descending Game Condition is exactly the Axiom of Foundation;
see Section~\ref{Sect:Impartial}.


\subsection{Adding and comparing games} \label{Ssec:Adding}
\cite[\S\S~7, 1]{ONAG}, \cite[\S~2]{WW}

Let us now introduce one of the most important concepts of the theory: the
{\em sum} of two games. Intuitively, we put two games next to each other
and allow each player to move in one of the two according to his choice,
leaving the other game unchanged; the next player can then decide
independently whether to move in the same game as her predecessor. The
negative of a game is the same game in which the allowed moves for both
players are interchanged (in games like chess, they simply switch colors).
The formal definitions are given below. Note that at this point it is
really necessary to require the DGC in its general form (rather than only
for alternating moves) in order to guarantee that the sum of two games is
again a game (which ends after a finite number of moves).

\begin{Definition}[Sum and Negative of Games]
     \label{DefSum} \sloppy \quad
     Let $G=\game{G^L, \dots}{G^R, \dots}$ and
     \hbox{$H=\game{H^L, \dots}{H^R, \dots}$}
     be two games. Then we define
     \begin{align*}
       G+H &:\equiv \game{G^L+H,G+H^L,\dots}{G^R+H,G+H^R,\dots} \,\,, \\
       -G  &:\equiv \game{-G^R,\dots}{-G^L,\dots}
           \qquad\text{and} \\
       G-H &:\equiv G+(-H) \,\,.
     \end{align*}
\end{Definition}
These are again recursive definitions. The
definition of $G+H$ requires knowing several sums of the form $G^L+H$
etc.\ which must be defined first. However, all these additions are easier
than $G+H$: recursive definitions work by induction without base, similarly
as Conway Induction (this time, for binary relations): the sum $G+H$ is
well-defined as soon as all options $G^L+H$ etc.\ are well-defined%
\footnote{
More formally, one could consider $G+H$ a formal pair of games and then
prove by Conway Induction that every such formal pair is in fact a game:
if all formal pairs $G^L+H$, $G+H^L$, $G^R+H$ and $G+H^R$ are games, then
clearly so is $G+H$. Similar remarks apply to the definition of
multiplication and elsewhere.
}.
To see how things get off the ground, note that the set of left options
of $G+H$ is
\begin{equation}
     \bigcup_{G^L}\left\{ G^L+H\right\} \cup
     \bigcup_{H^L}\left\{ G+H^L\right\}
     \label{EqAdditionLeft}
\end{equation}
where $G^L$ and $H^L$ run through the left options of $G$ and $H$. If $G$
and/or $H$ have no left options, then the corresponding unions are empty,
and there might be no left options of $G+H$ at all, or they might be all
of the form $G^L+H$ (or $G+H^L$). Therefore, $G+H$ and $-G$ are games.

As an example, $-1 \equiv \game{\;}{0}$ is really the negative of
$1 \equiv \game{0}{\;}$, justifying our notation. Also,
$\star + \star \equiv \game{\star}{\star}$, and the latter is easily
seen to be a zero game (whoever begins, loses), so $\star + \star = 0$.
The following properties justify the name `addition' for the operation
just defined.

\begin{Theorem} \label{ThmAddition}
     Addition is associative and commutative with $0\equiv\game{\;}{\;}$ as zero
     element. Moreover, all games $G$ and $H$ satisfy $-(-G)\equiv G$ and
     $-(G+H)\equiv (-G)+(-H)$.
\end{Theorem}
\begin{Proof}
     By (\ref{EqAdditionLeft}), the left (right) options of 
$G+\game{\;}{\;}$ are
     $G^L+\game{\;}{\;}$ ($G^R + \game{\;}{\;}$) only, so the claim
     `$G+\game{\;}{\;}\equiv G$' follows by Conway Induction.

     Commutativity uses induction too (in the second equality):
     \begin{align*}
       G+H
       &\equiv
       \game{G^L+H,G+H^L,\dots}{G^R+H,G+H^R,\dots} \\
       &\equiv
       \game{H+G^L,H^L+G,\dots}{H+G^R,H^R+G,\dots}
       \equiv
       H+G
       \,\,.
     \end{align*}

     Associativity works similarly (we write only left options):
     \begin{align*}
       (G+H)+K
       &\equiv
       \game{(G+H)^L+K,(G+H)+K^L,\dots}{\dots} \\
       &\equiv
       \game{(G^L+H)+K,(G+H^L)+K,(G+H)+K^L,\dots}{\dots} \\
       &\equiv
       \game{G^L+(H+K),G+(H^L+K),G+(H+K^L),\dots}{\dots} \\
       &\equiv
       \game{G^L+(H+K),G+(H+K)^L,\dots}{\dots}
       \equiv
       G+(H+K)
       \,\,.
     \end{align*}

     Moreover, omitting dots from now on,
     \[
       -(-G)
       \equiv
       -\game{-G^R}{-G^L}
       \equiv
       \game{-(-G^L)}{-(-G^R)}
       \equiv
       \game{G^L}{G^R}
       \equiv G
     \]
     where again induction was used in the third equality.
     Finally,
     \begin{align*}
       -(G+H)
       &\equiv
       -\game{G^L+H,G+H^L}{G^R+H,G+H^R} \\
       &\equiv
       \game{-(G^R+H),-(G+H^R)}{-(G^L+H),-(G+H^L)} \\
       &\equiv
       \game{(-G^R)+(-H),(-G)+(-H^R)}{(-G^L)+(-H),(-G)+(-H^L)} \\
       &\equiv
       \game{(-G)^L+(-H),(-G)+(-H)^L}{(-G)^R+(-H),(-G)+(-H)^R} \\
       &\equiv
       (-G)+(-H)
     \end{align*}
     where $-G^R$ means $-(G^R)$, etc. The third line uses induction again.
\end{Proof}

Conway calls inductive proofs like the preceding ones `one-line proofs'
(even if they do not fit on a single line): resolve the definitions,
apply induction, and plug in the definitions again.

Note that
\[ G - H :\equiv G+(-H) \equiv
\game{G^L-H, \dots, G-H^R, \dots}{G^R-H, \dots, G-H^L, \dots} \,. \]
   From now on, we will omit the dots in games like this (as already done in
the previous proof).

As examples, consider the games $2:\equiv
1+1\equiv \game{0+1,1+0}{\;}\equiv \game{1}{\;}$,
$3:\equiv 2+1\equiv \game{1+1,2+0}{\;}\equiv \game{2}{\;}$,
$4\equiv\game{3}{\;}$ etc., as well as
$-2\equiv\game{\;}{-1}$ etc.

\begin{Definition}
     We will write $G = H$ if $G - H =0$, $G > H$ if $G - H > 0$,
     $G \fu H$ if $G - H \fu 0$, etc.
\end{Definition}

It is obvious from the definition and the preceding result
that these binary relations extend
the unary relations $G = 0$ etc.\ defined earlier.

\begin{Theorem} \label{ThmFundamentalInequality}
     Every game $G$ satisfies $G=G$ or equivalently $G-G=0$. Moreover,
     $G^L \lf G$ for all left options $G^L$ and $G \lf G^R$ for all right
     options $G^R$ of $G$.
\end{Theorem}
\begin{Proof}
     By induction, we may suppose that $G^L-G^L\ge 0$ and $G^R-G^R\le 0$ for
     all left and right options of $G$. By definition, we have $G-G^R\ge 0$
     unless there is a right option $(G-G^R)^R\le 0$, and indeed such an
     option is $G^R-G^R\le 0$. Therefore, $G-G^R\lf 0$ or $G\lf G^R$. Similarly
     $G^L\lf G$.

     Now $G-G\equiv\game{G^L-G,G-G^R}{G^R-G,G-G^L}\ge 0$ unless any right option
     $(G-G)^R\le 0$; but we just showed that the right options are $G^R-G\gf 0$
     and $G-G^L\gf 0$, so indeed $G-G\ge 0$ and similarly $G-G\le 0$, hence
     $G-G=0$ and $G=G$.
\end{Proof}

The equality $G-G=0$ means that in the sum of any game with its negative,
the second player has a winning strategy: indeed, if the first player
makes any move in $G$, then the second player has the same move in $-G$
available and can copy the first move; the same holds if the first player
moves in $-G$ because $-(-G)\equiv G$. Therefore, the second player can
never run out of moves before the first does, so the Normal Play
Convention awards the win to the second player. This is sometimes
paraphrased like this: when playing against a Grand Master simultaneously
two games of chess, one with white and one with black, then you can
force at least one win (if draws are not permitted, as in our theory)!
In Mis\`ere Play, we would not have the fundamental equality $G=G$.

The following results show that the ordering of games is
compatible with addition.

\begin{Lemma} \label{LemTransitivity} \strut
     \begin{enumerate}
       \item If $G \ge 0$ and $H \ge 0$, then $G+H \ge 0$.
       \item If $G \ge 0$ and $H \gf 0$, then $G+H \gf 0$.
     \end{enumerate}
\end{Lemma}

Note that $G \gf 0$ and $H \gf 0$ implies nothing about $G+H$: the sum of
two fuzzy games can be in any outcome class. (Find examples!)

\begin{Proof}
     We prove both statements simultaneously using Conway Induction (with the
     binary relation $P(G,H)$: `for the pair of games $G$ and $H$, the
     statement of the Lemma holds'). The following proof can easily be
     rephrased in the spirit of `Left has a winning move unless\dots'.

     (1) $G\ge 0$ and $H\ge 0$ mean there are no $G^R\le 0$ and no $H^R\le 0$,
     so all $G^R\gf 0$ and all $H^R\gf 0$. By the inductive hypothesis, all
     $H+G^R\gf 0$ and $G+H^R\gf 0$, so $G+H$ has no right options $(G+H)^R\le
     0$ and thus $G+H\ge 0$.

     (2) Similarly, $H\gf 0$ means there is an $H^L\ge 0$. By the inductive
     hypothesis, $G+H^L\ge 0$, so $G+H$ has a left option $G+H^L\ge 0$ and thus
     $G+H\gf 0$.
\end{Proof}

\begin{Theorem} \label{ThmAddZero}
     The addition of a zero game never changes the outcome: if $G=0$, then
     $H>0$ or $H<0$ or $H=0$ or $H\fu 0$ iff $G+H>0$, $G+H<0$, $G+H=0$ or
     $G+H\fu 0$, respectively.
\end{Theorem}
\begin{Proof}
     If $H\ge 0$ or $H\le 0$, then $G+H\ge 0$ or $G+H\le 0$ by
     Lemma~\ref{LemTransitivity}, and similarly if $H\gf 0$ or $H\lf 0$, then
     $G+H\gf 0$ or $G+H\lf 0$. Since $H = 0$ is equivalent to $H \ge 0$
     and $H \le 0$, $H > 0$ is equivalent to $H \ge 0$ and $H \gf 0$, etc.,
     the `only if' direction follows. The `if'
     direction then follows from the fact that every game is in exactly one
     outcome class.
\end{Proof}

\begin{Corollary} \label{CorSameOutcome}
     Equal games are in the same outcome classes: if $G=H$, then $G>0$ iff
     $H>0$ etc.
\end{Corollary}
\begin{Proof}
     Consider $G+(H-H) \equiv H+(G-H)$, which by Theorem~\ref{ThmAddZero} has
     the same outcome class as $G$ and~$H$.
\end{Proof}

\begin{Corollary} \label{CorAddOrder}
     Addition respects the order: for any triple of games, $G>H$ is equivalent
     to $G+K>H+K$, etc.
\end{Corollary}
\begin{Proof}
     $G+K>H+K \Longleftrightarrow (G-H)+(K-K)>0 \Longleftrightarrow
     G-H>0 \Longleftrightarrow G>H$.
\end{Proof}

\begin{Theorem} \label{ThmOrder}
     The relation $\ge$ is reflexive, antisymmetric and transitive, and
     equality $=$ is an equivalence relation.
\end{Theorem}
\begin{Proof}
     Reflexivity of $\ge$ and $=$ is Theorem~\ref{ThmFundamentalInequality},
     and antisymmetry of $\ge$ and symmetry of $=$ are defined. Transitivity
     of $\ge$ and thus of $=$ follows like this:
     $G\geq H$ and $H\geq K$ implies $G-H\ge 0$ and $H-K\ge 0$, hence
     $G-K+(H-H)\ge 0$ by Lemma~\ref{LemTransitivity}. By
     Theorem~\ref{ThmAddZero}, this implies $G-K\ge 0$ and $G\ge K$.
\end{Proof}

\begin{Theorem} \label{ThmGroup}
     The equivalence classes formed by equal games form an additive abelian
     GROUP in which the zero element is represented by any game $G=0$.
\end{Theorem}
\begin{Proof}
     First we have to observe that addition and negation are compatible with
     respect to the equivalence relation: if $G=G'$ and $H=H'$ then $G-G'=0$
     and $H-H'=0$, hence $(G+H)-(G'+H')\equiv(G-G')+(H-H')=0$ by
     Lemma~\ref{LemTransitivity} and $G+H=G'+H'$ as needed. Easier yet, $G=G'$
     implies $0=G-G'\equiv -(-G)+(-G')\equiv(-G')-(-G)$, hence $-G'=-G$.

     For every game $G$, the game $-G$ represents the inverse equivalence class
     by Theorem~\ref{ThmFundamentalInequality}. Finally, addition is
     associative and commutative by Theorem~\ref{ThmAddition}.
\end{Proof}

It is all well to define equivalence classes of games, but their
significance sits in the fact that replacing a game by an equivalent one
never changes the outcome, even when this happens for games that are
themselves parts of other games.

\begin{Theorem}[Equal Games]
\label{ThmEqualGames}
     If $H=H'$, then $G+H=G+H'$ for all games~$G$.
     If $G=\game{G^{L_1},G^{L_2},\dots}{G^{R_1},G^{R_2},\dots}$
     and $H=\game{H^{L_1},H^{L_2},\dots}{H^{R_1},H^{R_2},\dots}$ are two games
     such that $G^{L_i}=H^{L_i}$ and $G^{R_i}=H^{R_i}$ for all left and right
     options, then $G=H$:
     replacing any option by an equivalent one (or any set of options by
     equivalent options) yields an equivalent game.
\end{Theorem}
\begin{Proof}
     The first part is self-proving: $(G+H)-(G+H')=(G-G)+(H-H')=0$. The second
     part is similar, but easier to write in words: in $G-H$, Left might move
     in $G$ to some $G^L-H$ or in $H$ to some $G-H^R$, and Right's answer will
     be either in $H$ to a $G^L-H^{L'}=0$ (with $H^{L'}$ chosen so that
     $H^{L'}=G^L$), or Right answers in $G$ to a $G^{R'}-H^R=0$. The
     situation is analogous if Right starts.
\end{Proof}


\subsection{Simplifying games}
\label{Ssect:Simplifying}
\cite[\S~3]{WW}, \cite[\S~10]{ONAG}

Since equality of games is a defined equivalence
relation, there are many ways of writing down a game that has a
certain value (i.e., lies in a certain equivalence class). Some of
these will be simpler than others, and there may even be a simplest
or canonical form of a game. In this section, we show how one can
simplify games and that simplest forms exist for an interesting
class of games.

\begin{Definition}[Gift Horse]
     Let $G$ and~$H$ be games. If $H \lf G$, then $H$ is a {\em left
     gift horse} for~$G$; if $H \gf G$, then $H$ is a {\em right gift
     horse} for~$G$
\end{Definition}

\begin{Lemma}[Gift Horse Principle] \label{LemmaGiftHorse}
     If $H_L, \dots$ are left gift horses and $H_R, \dots$ are
     right gift horses for~$G=\game{G^L,\dots}{G^R,\dots}$, then
     \[ \game{H_L, \dots, G^L, \dots}{H_R, \dots, G^R, \dots} = G \,. \]
     (Here, $\{H_L, \dots\}$ and $\{H_R, \dots\}$ can be arbitrary sets
     of games.)
\end{Lemma}
\begin{Proof}
     Let
     $G' \equiv \game{H_L, \dots, G^L, \dots}{H_R, \dots, G^R, \dots}$.
     Then $G' - G \ge 0$, since the right options are $G^R - G \gf 0$
     (by Theorem~\ref{ThmFundamentalInequality}), $H_R - G \gf 0$
     (by assumption),
     and $G' - G^L$, which has the left option $G^L - G^L = 0$,
     so $G' - G^L \gf 0$. In the same way, we see that $G' - G \le 0$,
     and it follows that $G' = G$.
\end{Proof}

This `Gift Horse Principle' tells us how to offer extra options to a
player without changing the value of a game (since no player really wants
to move to these options), so we know how to make games more complicated.
Now we want to see how we can {\em remove} options and thereby make
a game simpler. Intuitively, an option that is no better than another
option can as well be left out, since a reasonable player will never
use it. This is formalized in the following definition and lemma.

\begin{Definition}[Dominated Option]
     Let $G$ be a game. A left option~$G^L$ is {\em dominated} by another left
     option~$G^{L'}$ if $G^L \le G^{L'}$. Similarly, a right option~$G^R$ is
     {\em dominated} by another right option~$G^{R'}$ if $G^R \ge G^{R'}$.
\end{Definition}

\begin{Lemma}[Deleting Dominated Options] \label{LemmaDominatedOptions}
     Let $G$ be a game with fixed left and right options $G^L$ and $G^R$. Then
     the value of $G$ remains unchanged if some or all left options which are
     dominated by $G^L$ are removed, and similarly if some or all right
     options which are dominated by $G^R$ are removed.
\end{Lemma}
\begin{Proof}
     Let $G'$ be
     the game obtained from~$G$ by removing all or some left options
     that are dominated by~$G^L$ (but keeping~$G^L$ itself). Then all
     the deleted options are left gift horses for~$G'$, since for such
     an option~$H$, we have $H \le G^L \lf G'$. We can therefore add all
     these options to~$G'$, thereby obtaining~$G$, without changing the value.
     The same argument works for dominated right options.
\end{Proof}

As simple examples, we have $2\equiv\game{1}{\;}=\game{0,1}{\;}$,
$3\equiv\game{2}{\;}=\game{0,1,2}{\;}$ etc., so we recover von Neumann's
definition of natural numbers. Another example would be $\game{0,1}{2,3} =
\game{1}{2}$. Note that it is possible that all options are dominated,
but this does not mean that all options can be removed: as an example,
consider
$\omega :\equiv \game{0, 1, 2, \dots}{\;}$.

There is another way of simplifying a game that does not work by
removing options, but by introducing shortcuts. The idea is as follows.
Suppose Left has a move~$G^L$ which Right can counter to some fixed
$G^{LR} \le G$, a position at least as good for Right as $G$ was.
The claim is that replacing the single option $G^L$ by all left options
of $G^{LR}$ does not change the value of $G$: this does not hurt Left (if
Left wants to move to $G^L$, then he must expect the answer $G^{LR}$ and
then has all left options of $G^{LR}$ available); on the other hand, it
does not help Left if $G$ is replaced by $G^{LR}\le G$. Precise
statements are like this.

\begin{Definition}[Reversible Option]
     Let $G$ be a game. A left option~$G^L$ is called {\em reversible
     (through~$G^{LR}$)} if $G^L$ has a right option $G^{LR} \le G$.
     Similarly, a right option~$G^R$ is called {\em reversible
(through~$G^{RL}$)}
     if $G^R$ has a left option $G^{RL} \ge G$.
\end{Definition}

\begin{Lemma}[Bypassing Reversible Options] \label{LemmaReversible}
     If $G$ has a left option~$H$ that is reversible through~$K = H^R$,
     then
     \[ G = \game{H, G^L, \dots}{G^R, \dots}
          = \game{K^L, \dots, G^L, \dots}{G^R, \dots}
     \]
(here, $G^L$ runs through all left options of $G$ other than $H$).
     In words: the value of~$G$ is unchanged when we replace the reversible
     left option~$H$ by all the left options of~$K$. A similar statement
     holds for right options.
\end{Lemma}
\begin{Proof}
     Let $G' = \game{K^L, \dots, G^L, \dots}{G^R, \dots}$ and
     $G'' = \game{H, K^L, \dots, G^L, \dots}{G^R, \dots}$.
     We claim that $H$ is a left gift horse for~$G'$. This can be seen as
     follows. First, for all~$K^L$ we have $K^L \lf G'$, since $K^L$ is a
     left option of~$G'$. Also, $K \le G \lf G^R$, so $K \lf G^R$ for all~$G^R$.
     These statements together imply that $K \le G'$. Since $K$ is a right
     option of~$H$, this in turn says that $H \lf G'$, as was to be shown.
     By the Gift Horse Principle, we now have $G' = G''$. On the other
     hand, $K^L \lf K \le G$, so all the $K^L$ are left gift horses for~$G$,
     whence $G = G'' = G'$.
\end{Proof}

One aspect of reversible options might be surprising: if $G^L$ is
reversible through $G^{LR}$, this means that Left may bypass the move to
$G^L$ and Right's answer to $G^{LR}$ and move directly to some left
option of $G^{LR}$; but what if there was another right option $G^{LR'}$
which Right might prefer over $G^{LR}$: is Right deprived of her better
move? For the answer, notice that $\game{\;}{1}=\game{\;}{100}=0$:
although Right might prefer that her only move was $1$ rather than $100$,
the first player to move will always lose, which is all that counts.
Similarly, depriving Right of her better answer $G^{LR'}$ would make a
difference only if there was a game $S$ such that $G+S\le 0$ but
$G^{LR}+S\gf 0$ (our interest is in the case that Left starts:
these conditions mean that Left cannot win in $G+S$, but he can when
jumping directly to $G^{LR}$); however, the first condition and the
hypothesis imply $G^{LR}+S\le G+S\le 0$, contradicting the second
condition.

Given the simplifications of games described above, the question arises
whether there is a simplest form of a game: a form that cannot be further
simplified by removing dominated options and bypassing reversible
options. The example  $\omega = \game{0, 1, 2, \dots}{\;}$ shows that this
is not the case in general. But such a simplest form exists if we impose
a natural finiteness condition which is satisfied by most real-life games.

\begin{Definition}[Short]
      A game~$G$ is called \emph{short} if it has only finitely many positions.
\end{Definition}

\begin{Theorem}[Normal Form] \label{ThmShortNormalForm}
     In each equivalence class of short games, there is a unique game
     that has no dominated or reversible positions.
\end{Theorem}
\begin{Proof}
     Since both ways of simplifying games reduce the number of positions,
     we eventually reach a game that cannot be simplified further. This
     proves existence.

     To prove uniqueness, we assume that $G$ and~$H$ are two equal (short)
     games both without dominated and reversible positions. We have to show
     that $G \equiv H$. Let $G^L$ be some left option of~$G$. Since
     $G^L \lf G = H$, there must be a right option~$G^{LR} \le H$
     or a left option~$H^L$ such that $G^L \le H^L$. The first is
     impossible since $G^L$ is not reversible. Similarly, there is some~$G^{L'}$
     such that $H^L \le G^{L'}$, so $G^L \le G^{L'}$. But there are no
     dominated options either, so $G^L = H^L = G^{L'}$. By induction,
     $G^L \equiv H^L$. In that way, we see that $G$ and~$H$ have the same
     set of (identical) left options, and the same is true for the right
     options.
\end{Proof}


\section{The FIELD of Numbers}
\label{Sect:Field}


\subsection{What is a number?} \label{Ssect:Number}
\cite[\S\S~0,1]{ONAG}, \cite[\S~2]{WW}

We already have encountered games like $0$, $1$, $-1$, $2$ that we
have denoted by numbers and that behave like numbers. In particular,
they measure which player has got how many free moves left and therefore
are easy to compare. We now want to extend this to a class of games that
is as large as possible (and whose elements are to be called {\em numbers}).

The guiding idea is that numbers should be totally ordered, i.e.\ no two
numbers should ever be fuzzy to each other.
Recall that by Thm.~\ref{ThmFundamentalInequality} we always have
$G^L \lf G$ and $G \lf G^R$. If $G$, $G^L$ and~$G^R$ are to be numbers,
this forces $G^L < G < G^R$, so we must at least require that $G^L < G^R$.
In order for numbers to be preserved under playing, we need to require
that all options of numbers are numbers. This leads to the following
definition.

\begin{Definition}[Number] \label{DefNumber}
     A game $x=\game{x^L, \dots}{x^R, \dots}$ is called a {\em number}
if all left
     and right options $x^L$ and $x^R$ are numbers and satisfy $x^L<x^R$.
\end{Definition}

As it turns out, this simple definition leads not only to a totally ordered
additive subGROUP of games but even to an algebraically closed FIELD which
simultaneously contains the real and ordinal numbers!

We will use lowercase letters $x, y, z, \dots$ to denote numbers.
The simplest numbers are $0$, $1$ and $-1$. A slightly more interesting
number is $\tfrac{1}{2} :\equiv \game{0}{1}$ (one checks easily that
$\tfrac{1}{2} + \tfrac{1}{2} = 1$, justifying the name). There are
also `infinite numbers' like $\omega = \game{0, 1, 2, \dots}{\;}$.

\begin{Lemma} \label{LemFundamentalInequalityNumbers}
     Every number $x=\game{x^L, \dots}{x^R, \dots}$ satisfies $x^L<x<x^R$.
\end{Lemma}
\begin{Proof}
     The left options of $x^L-x$ are of the form $x^L-x^R$ or $x^{LL}-x$.
     Since $x$ is a number, we have $x^L-x^R<0$. We use the inductive
     hypothesis $x^{LL}<x^L$ and $x^L-x\lf 0$ from
     Theorem~\ref{ThmFundamentalInequality}. Therefore,
     Lemma~\ref{LemTransitivity} implies $x^{LL}-x=(x^{LL}-x^L)+(x^L-x)\lf 0$.

     If $x^L-x\gf 0$ was true, we would need some $(x^L-x)^L\ge 0$,
     which we just excluded. Hence $x^L\le x$ for all left options $x^L$ of
     $x$, and similarly $x\le x^R$ for all right options $x^R$. The claim
     now follows because $x^L\lf x\lf x^R$ from
     Theorem~\ref{ThmFundamentalInequality}.
\end{Proof}

\begin{Theorem} \label{ThmAdditionNumbers}
     If $x$ and $y$ are numbers, then $x+y$ and $-x$ are numbers, so
     (equivalence classes of) numbers form an abelian subGROUP of games.
\end{Theorem}
\begin{Proof}
     Since $-x=\game{-x^R, \dots}{-x^L, \dots}$,
     we have $(-x)^L=-x^R<-x^L=(-x)^R$,
     so the options of $-x$ are ordered as required. Conway Induction now shows
     that $-x$ is a number.

     In $x+y=\game{x^L+y,x+y^L,\dots}{x^R+y,x+y^R,\dots}$, we have the
inequalities
     $x^L+y<x^R+y$ and $x+y^L<x+y^R$ by Corollary~\ref{CorAddOrder}. By
     Lemma~\ref{LemFundamentalInequalityNumbers}, we also
     have $x^L+y<x+y<x+y^R$ and $x+y^L<x+y<x^R+y$, so $x+y$ is a number as soon
     as all its options are, and Conway Induction applies.
\end{Proof}

\begin{Theorem} \label{ThmNumbersOrdered}
     Numbers are totally ordered: every pair of numbers $x$ and $y$ satisfies
     exactly one of $x<y$, $x>y$, or $x=y$.
\end{Theorem}
\begin{Proof}
     Suppose there was a number $z\fu 0$. This would imply the existence of
     options $z^L\ge 0\ge z^R$, which is excluded by definition: numbers are
     never fuzzy.

     Now if there were two numbers $x\fu y$, then $x-y$ would be a number by
     Theorem~\ref{ThmAdditionNumbers} and $x-y\fu 0$, but this is impossible,
     as we have just shown.
\end{Proof}


\subsection{Short numbers and real numbers}
\label{Ssect:Short}
\cite[\S~2]{ONAG}, \cite[\S~2]{WW}

A short number is simply a number that is a short game, i.e.\ a game
with only finitely many positions. In particular, it then has only
finitely many options, and since numbers are totally ordered, we can
eliminate dominated options so as to leave at most one left, resp.\
right option.

By the definition of negation, addition and multiplication (see below
in Section~\ref{Ssect:Multiplication}), it is
easily seen that the set~(!) of (equivalence classes of) short numbers
forms a unitary ring.

\begin{Theorem}
     The ring of short numbers is (isomorphic to) the ring~$\Z[\frac{1}{2}]$
     of dyadic fractions.
\end{Theorem}
\begin{Proof}
     We have already seen that $\game{0}{1} = \frac{1}{2}$, therefore
     $\Z[\frac{1}{2}]$ is contained in the ring of short numbers.
     For the converse, see~\cite[Theorem~12]{ONAG}.
     The main step in proving the converse is to show that
     \[ \left\{\frac{m}{2^n}\biggm|\frac{m+1}{2^n}\right\}
          = \frac{2m+1}{2^{n+1}}
     \]
     for integers $m$ and natural numbers $n$.
\end{Proof}

Let $S$ denote the ring of short numbers (or dyadic fractions). We can
represent every element~$x$ of~$S$ in the form
\[ x = \game{y \in S : y < x}{y \in S : y > x} \,, \]
where both sets of options are nonempty.
In fact, the set of all numbers satisfying this property
is exactly the field~$\mathbb R$ of real numbers: we are taking
Dedekind sections in the ring~$S$. (More precisely, this is the most
natural model of the real numbers within our Conway numbers:
it is the only one where all real numbers have all their options in $S$.
There are other embeddings that are obtained by choosing a transcendence basis
of~$\mathbb R$ over $\mathbb Q$ and then changing the images
of this basis by some infinitesimal amounts.)
For some more discussion,
see~\cite[Chapter~2]{ONAG}.


\subsection{Multiplication of numbers}
\label{Ssect:Multiplication}

In order to turn numbers into a FIELD, we need a multiplication.

\begin{Definition}[Multiplication] \label{DefMultiplication}
     Given two numbers $x=\game{x^L, \dots}{x^R, \dots}$ and
     \hbox{$y=\game{y^L, \dots}{y^R, \dots}$,} we define the product
     \begin{eqnarray*}
       x\cdot y
       &:=&
       \gameL{x^L\cdot y+x\cdot y^L-x^L\cdot y^L,\;
              x^R\cdot y+x\cdot y^R-x^R\cdot y^R,\dots} \\
       && \quad \gameR{x^L\cdot y+x\cdot y^R-x^L\cdot y^R,\;
                       x^R\cdot y+x\cdot y^L-x^R\cdot y^L,\dots}
     \end{eqnarray*}
     As with addition in (\ref{EqAdditionLeft}), the left and right options
     are all terms as in the definition that can be formed with left and right
     options of $x$ and $y$. More precisely, the left options are indexed
     by pairs $(x^L, y^L)$ and pairs $(x^R, y^R)$, and similarly for the
     right options.
     In particular, this means that if $x$ has no left options (say), then
     $x \cdot y$ will not have left and right options of the first type shown.
     We will usually omit the dot and write $xy$ for $x\cdot y$.
\end{Definition}

While this definition might look complicated, it really is not. It is
motivated by $x^L<x<x^R$ and $y^L<y<y^R$, so we want multiplication to
satisfy $(x-x^L)(y-y^L)>0$, hence $xy>x^Ly+xy^L-x^Ly^L$, which motivates
the first type of left options. The other three types are obtained in a
similar way.

One might try the simpler definition
$xy=\game{x^Ly,xy^L}{x^Ry,xy^R}$ for multiplication, motivated by
$x^L<x<x^R$ and $y^L<y<y^R$. But the inequalities would be wrong for
negative numbers. In fact, this would be just a different notation for
addition!

Recall the two special numbers $0\equiv\game{\;}{\;}$ and
$1\equiv\game{0}{\;}$.

\begin{Theorem} \label{ThmMultiplication1}
     For all numbers $x,y,z$, we have the identities
     \[
       0\cdot x\equiv 0\,,
       \quad
       1\cdot x\equiv x\,,
       \quad
       x y\equiv y x\,,
       \quad
       (-x) y \equiv y(-x) \equiv -xy
     \]
     and the equalities
     \[
       (x+y) z = xz + yz\,, \qquad (xy)z=x(yz)
       \,\,.
     \]
\end{Theorem}
\begin{Proof}
     The proofs are routine `1-line-proofs'; the last two are a bit
     lengthy and can be found in \cite[Theorem~7]{ONAG}.
\end{Proof}

The reason why multiplication is defined only for numbers, not for
arbitrary games, is that there are games $G,G',H$ with $G=G'$ but
$GH\neq G'H$. For example, we have $\game{1}{\;} = \game{0,1}{\;}$,
but $\game{1}{\;} \cdot \star = \game{\star}{\star} = 0$, whereas
$\game{0,1}{\;} \cdot \star = \game{0, \star}{0, \star} \fu 0$. (Note
that we have $0 \cdot G = 0$ and $1 \cdot G = G$ for any game~$G$.
Furthermore, since games form an abelian GROUP, we always have integral
multiples of arbitrary games.)

The following theorem shows that our multiplication behaves as expected.
Its proof is the most complicated inductive proof in this
paper. It required quite some work to produce a concise version of the
argument in the proof, even with Conway's \cite[Theorem~8]{ONAG} at hand.
The main difficulty is to organize a simultaneous induction for three
different statements with different arguments.

\begin{Theorem} \label{ThmMultiplication2} \strut
     \begin{enumerate}
       \item \label{TM2i1}
         If $x$ and $y$ are numbers, then so is $xy$.
       \item \label{TM2i2}
         If $x_1$, $x_2$, $y$ are numbers such that $x_1 = x_2$,
         then $x_1 y = x_2 y$.
       \item \label{TM2i3}
         If $x_1$, $x_2$, $y_1$, $y_2$ are numbers such that $x_1 < x_2$ and
         $y_1 < y_2$, \\ \quad then
         $x_1 y_2 + x_2 y_1 < x_1 y_1 + x_2 y_2$. \\
         In particular, if $y > 0$, then $x_1 < x_2$ implies $x_1 y < x_2 y$.
     \end{enumerate}
\end{Theorem}
\begin{Proof}
     We will prove most of the statements simultaneously using Conway
Induction.
     More precisely, let $P_1(x, y)$, $P_2(x_1, x_2, y)$ and
     $P_3(x_1, x_2, y_1, y_2)$ stand for the statements above.
     For technical reasons, we also introduce the statement
     $P_4(x, y_1, y_2)$:
     \begin{enumerate} \addtocounter{enumi}{3}
       \item \label{TM2i4}
        {\em  If $x$, $y_1$, $y_2$ are numbers with $y_1<y_2$, then
        $P_3(x^L, x, y_1, y_2)$ and $P_3(x, x^R, y_1, y_2)$ hold for all
        options $x^L$, $x^R$ of~$x$.}
     \end{enumerate}

     We begin by proving $P_1$, $P_2$ and~$P_4$ simultaneously.
     For this part, we assume that all occurring numbers are
     positions in a fixed number~$z$ (we can take
     $z = \game{x,x_1, x_2,y,y_1,y_2}{\;}$). Then for a position $z'$ in~$z$,
     we define $n(z')$ to be the distance between $z'$ and the root~$z$
     in the rooted tree representing the game~$z$: that is the number of
    moves needed to reach $z'$ from the starting position $z$ (a non-negative
    integer). Formally, we set
     $n(z) = 0$ and, if $z'$ is a position in~$z$, $n({z'}^L) = n(z') + 1$,
     $n({z'}^R) = n(z') + 1$.

     For each statement $P_1$, $P_2$, $P_4$, we measure its `depth'
     by a pair of natural numbers $(r, s)$ (where $s = \infty$ is allowed)
     as follows.
     \begin{itemize}
       \item The depth of $P_1(x, y)$ is $(n(x) + n(y), \infty)$.
       \item The depth of $P_2(x_1, x_2, y)$
             is $(\min\{n(x_1), n(x_2)\} + n(y), \max\{n(x_1), n(x_2)\})$.
       \item The depth of $P_4(x, y_1, y_2)$
             is $(n(x) + \min\{n(y_1), n(y_2)\}, \max\{n(y_1), n(y_2)\})$.
     \end{itemize}
     The inductive argument consists in showing that each statement follows
     from statements that have greater depths (in the lexicographic
     ordering) and only involve positions of the games occurring in the
     statement under consideration. If the statement was false, it would
     follow that there was an infinite chain of positions $z'$ in~$z$ of
     unbounded depth~$n(z')$, each a position of its predecessor; this would
     contradict the Descending Game Condition.

     Properties of addition will be used without explicit mention.

     (1) We begin with $P_1(x, y)$. We may assume that all terms
     $x^L y$, $x y^L$, $x^L y^L$ etc.\ are numbers, using $P_1(x^L, y)$ etc.;
     therefore
     all options of $x y$ are numbers. It remains to show that the
     left options are smaller than the right options. There are four
     inequalities to show; we treat one of them in detail (the others being
     analogous). We show that
     \[ x^{L_1} y + x y^L - x^{L_1} y^L < x^{L_2} y + x y^R - x^{L_2} y^R \,. \]
     There are three cases. First suppose that $x^{L_1} = x^{L_2}$.
     Then using $P_2(x^{L_1}, x^{L_2}, y)$ and
     $P_2(x^{L_1}, x^{L_2}, y^R)$, the statement is equivalent to
     $P_3(x^{L_1}, x, y^L, y^R)$, which is in turn a special case
     of~$P_4(x, y^L, y^R)$.

     Now suppose that $x^{L_1} < x^{L_2}$. Then we use
     $P_3(y^L, y, x^{L_1}, x^{L_2})$, which follows
     from $P_4(y, x^{L_1}, x^{L_2})$,
     and $P_3(x^{L_2}, x, y^L, y^R)$,
     which follows from $P_4(x, y^L, y^R)$, to get
     \[
       x^{L_1} y + x y^L - x^{L_1} y^L < x^{L_2} y + x y^L - x^{L_2} y^L
                                       < x^{L_2} y + x y^R - x^{L_2} y^R \,.
     \]

     Similarly, if $x^{L_1} > x^{L_2}$, we use $P_4(x, y^L, y^R)$
     and~$P_4(y, x^{L_2}, x^{L_1})$ to get
     \[
       x^{L_1} y + x y^L - x^{L_1} y^L < x^{L_1} y + x y^R - x^{L_1} y^R
                                       < x^{L_2} y + x y^R - x^{L_2} y^R \,.
     \]

     (2) For~$P_2(x_1, x_2, y)$, note that $z_1 = z_2$ if $z_1^L < z_2 <
     z_1^R$ and $z_2^L < z_1 < z_2^R$ for all relevant options. So
     we have to show a number of statements of the
     type $(x_1 y)^L < x_2 y$ or $ x_2 y < (x_1 y)^R$.
                        We carry this out for the left option
     $(x_1 y)^L = x_1^L y + x_1 y^L - x_1^L y^L$; the other possible
     cases are done in the same way. Statement $P_2(x_1, x_2, y^L)$ gives
     $x_1 y^L = x_2 y^L$ and $P_4(y, x_1^L, x_2)$ gives
     $x_1^L y + x_2 y^L < x_1^L y^L + x_2 y$, which together imply
     $(x_1 y)^L < x_2 y$.

     (3) We now consider $P_4(x, y_1, y_2)$. Since $y_1 < y_2$, there is
     some $y_1^R$ such that $y_1 < y_1^R \le y_2$, or there is some~$y_2^L$
     such that $y_1 \le y_2^L < y_2$. We consider the first case; the second
     one is analogous. First note that from $P_1(x, y_1)$, we get
     $(x y_1)^L < x y_1 < (x y_1)^R$
     for all left and right options of~$x y_1$. We therefore obtain the
     inequalities
     \begin{equation}
        x y_1 + x^L y_1^R < x^L y_1 + x y_1^R \And
        x^R y_1 + x y_1^R < x y_1 + x^R y_1^R \,.
          \label{Eq:NumbersProof1}
     \end{equation}
     Now if $y_1^R = y_2$, then by $P_2(y_1^R, y_2, x)$, $P_2(y_1^R, y_2, x^L)$
     and $P_2(y_1^R, y_2, x^R)$ we are done.
     Otherwise, $y_1^R<y_2$ and $P_4(x, y_1^R, y_2)$ says
     \begin{equation}
        x^L y_2 + x y_1^R < x^L y_1^R + x y_2 \And
        x y_2 + x^R y_1^R < x y_1^R + x^R y_2 \,.
          \label{Eq:NumbersProof2}
     \end{equation}
     Adding the left resp.\ right inequalities in (\ref{Eq:NumbersProof1})
     and (\ref{Eq:NumbersProof1}) and canceling like terms proves the claim.

     This shows that every statement $P_1(x,y)$, $P_2(x_1,x_2,y)$, and
     $P_4(x,y_1,y_2)$ follows from similar statements which use the same
     arguments or some of their options, and it is easily verified that all
     used statements have greater depths. This proves $P_1$, $P_2$ and~$P_4$.

     (4) It remains to show~$P_3(x_1, x_2, y_1, y_2)$. This is done by
     Conway Induction in the normal way, using the statements we have
     already shown.

     Since $x_1 < x_2$, there is
     some $x_1^R \le x_2$ or some $x_2^L \ge x_1$. Assume the first possibility
     (the other one is treated in the same way). If $x_1^R = x_2$,
     then we apply $P_2$ to get $x_1^R y_1 = x_2 y_1$ and $x_1^R y_2 = x_2 y_2$.
     Using $P_4(x_1, y_1, y_2)$, we also have
     $x_1 y_2 + x_1^R y_1 < x_1 y_1 + x_1^R y_2$; together, these
     imply the desired conclusion. Finally, if $x_1^R < x_2$, then
     by induction ($P_3(x_1^R, x_2, y_1, y_2)$),
     we get $x_1^R y_2 + x_2 y_1 < x_1^R y_1 + x_2 y_2$ and, using
     $P_4(x_1, y_1, y_2)$ again,
     also $x_1 y_2 + x_1^R y_1 < x_1 y_1 + x_1^R y_2$.
     Adding them together and canceling like terms proves our claim.
\end{Proof}


\subsection{Division of numbers}
\label{Ssect:Division}

The definition of division is more complicated than that of addition or
multiplication --- necessarily so, since for example $3=\game{2}{\;}$ is a
very simple game with only finitely many positions, while
\(
\frac{1}{3} =
\game{\frac{1}{4}, \frac{5}{16}, \frac{21}{64}, \dots}%
        {\dots, \frac{11}{32}, \frac{3}{8}, \frac{1}{2}}
\)
has infinitely many positions which must all be `generated' somehow
from the positions of~$3$.

It suffices to find a multiplicative inverse for every $x>0$. It is
convenient to rewrite positive numbers as follows.

\begin{Lemma}
     For every number $x>0$, there is a number $y$ without negative options
     such that $y=x$.
\end{Lemma}
\begin{Proof}
     To achieve this, we simply add the left Gift Horse $0$ and then delete all
     negative left options, which are now dominated by $0$.
\end{Proof}

\begin{Theorem}
     For a number $x>0$ without negative options, define
     \[
       y=\gameBig{0, \frac{1 + (x^R-x) y^L}{x^R}, \frac{1 + (x^L-x) y^R}{x^L}}%
                 {\frac{1 + (x^L-x) y^L}{x^L},\frac{1 + (x^R-x) y^R}{x^R}}
     \]
     where all options $x^L\neq 0$ and $x^R$ of $x$ are used.
     Then $y$ is a number with $xy=1$.
\end{Theorem}
Note that this definition is recursive as always: in order to find
$y=1/x$, we need to know $1/x^L$ and $1/x^R$ first. However, this time we
also need to know left and right options of $y$. We view this really as
an algorithmic definition: initially, make $0$ a left option of $y$. Then,
for every left option $y^L$ generated so far, produce new left and right
options of $y$ with all $x^L$ and all $x^R$; similarly, for every right
option $y^R$ already generated, do the same. This step is then iterated
(countably) infinitely often. More precisely, we define a sequence
of pairs of sets of numbers $Y^L_n, Y^R_n$ recursively as follows.
\begin{align*}
      Y^L_0 &= \{0\}, \, \quad Y^R_0 = \emptyset \\
      Y^L_{n+1} &= Y^L_n
      \cup \bigcup_{x^R} \Bigl\{\frac{1 + (x^R-x) y^L}{x^R} : y^L \in
Y^L_n\Bigr\}
      \cup \bigcup_{x^L} \Bigl\{\frac{1 + (x^L-x) y^R}{x^L} : y^R \in
Y^R_n\Bigr\}
          \\
      Y^R_{n+1} &= Y^R_n
      \cup \bigcup_{x^L} \Bigl\{\frac{1 + (x^L-x) y^L}{x^L} : y^L \in
Y^L_n\Bigr\}
      \cup \bigcup_{x^R} \Bigl\{\frac{1 + (x^R-x) y^R}{x^R} : y^R \in
Y^R_n\Bigr\}
\end{align*}
Then
\[ y=\gameBig{\bigcup_{n\in\N} Y^L_n}{\bigcup_{n\in\N} Y^R_n} \,\,. \]

If $x$ is a real number (which always can be written with at most
countably many options), we generate in every step new sets of left and
right options, the suprema and infima of which converge to $1/x$. In this
case, we have a convergent (infinite) algorithm specifying a Cauchy
sequence of numbers. General numbers $x$ might need arbitrarily big {\em
sets} of options, so the option sets of $y$ can become arbitrarily big
too. However, the necessary number of iteration steps in the construction
of $y$ is still at most countable.

\begin{Proof}
     (Compare \cite[Theorem~10]{ONAG}.)
     By induction, all the options of~$y$ are numbers (there are of course
     two inductive processes involved, one with respect to~$x$, and the
     other with respect to~$n$ above). We now prove that we have
     $x y^L < 1 < x y^R$ for all left and right options of~$y$. This is done
     by induction on~$n$. The statement is obvious for $n = 0$. We give one
     of the four cases for the inductive step in detail: suppose
     $y^L = (1 + (x^R-x) y^{L'})/x^R$ with $y^{L'} \in Y^L_n$, then by induction
     we have $x y^{L'} < 1$, so (by Theorem~\ref{ThmMultiplication2})
     $(x^R - x) x y^{L'} < x^R - x$, which is equivalent to the claim.

     In order to prove that $y$ is a number, we must show
     that the left options are smaller than the right options.
     It is easy to see that the right options are positive: for
     $(1 + (x^L - x)y^L)/x^L$, this follows from $1 - xy^L > 0$;
     for $(1 + (x^R - x)y^R)/x^R$, it follows by induction ($y^R > 0$).
     There are then
     four more cases involving the two different kinds of generated
     left resp.\ right options. We look at two of them, the other two are
     dealt with analogously. So suppose we want to show that
     $(1 + (x^R-x) y^{L_1})/x^R < (1 + (x^L-x) y^{L_2})/x^L$. This is equivalent
     to
     \[ x^R (1 + (x^L-x) y^{L_2}) - x^L (1 + (x^R-x) y^{L_1}) > 0 \,. \]
     The left hand side of this equation can be written in each of the
     following two ways:
     \begin{gather*}
        (x^R - x^L)(1 - x y^{L_1}) + (y^{L_1} - y^{L_2}) x^R (x - x^L) \\
        \hspace{10ex}
         {} = (x^R - x^L)(1 - x y^{L_2}) + (y^{L_2} - y^{L_1}) x^L (x^R - x) \,,
     \end{gather*}
     showing that it is always positive. Now suppose we want to show that \\
     $(1 + (x^{R_1}-x) y^L)/x^{R_1} < (1 + (x^{R_2}-x) y^R)/x^{R_2}$.
     This is equivalent to
     \[ x^{R_1} (1 + (x^{R_2}-x) y^R) - x^{R_2} (1 + (x^{R_1}-x) y^L) > 0 \,. \]
     Again, the left hand side can be written as
     \begin{gather*}
        (y^R - y^L) x^{R_1} (x^{R_2} - x) + (x^{R_1} - x^{R_2}) (1 - x y^L) \\
        \hspace{10ex}
        {} = (y^R - y^L) x^{R_2} (x^{R_1} - x) + (x^{R_2} - x^{R_1}) (x y^R - 1)
        \,,
     \end{gather*}
     showing that it is always positive. Note that we are using the inductive
     result that the `earlier' options $y^R$ and $y^L$ satisfy
     $y^R > y^L$.

     Finally, to prove that $xy = 1$, we have to show that $(xy)^L < 1 < (xy)^R$
     (since $0 = 1^L < xy$ trivially). For example, take
     \begin{align*}
       (xy)^R &= x^R y + x y^L - x^R y^L \\
              &= 1 + x^R \left(y - \frac{1 + (x^R - x) y^L}{x^R}\right) \\
              &= 1 + x^R (y - y^{L'}) > 1 \,.
     \end{align*}
\end{Proof}

\begin{Corollary}
\label{CorField}
     The (equivalence classes of) numbers form a totally ordered FIELD.
\end{Corollary}
In fact, this field is real algebraically closed. This is shown in
\cite[Chapter~4]{ONAG}.


\section{Ordinal numbers}
\label{Sect:Ordinal}


\subsection{Ordinal Numbers}
\label{Ssect:Ordinal}
\cite[\S~2]{ONAG}

\begin{Definition}[Ordinal Number]
      A game $G$ is an \emph{ordinal number} if it has no right options
      and all of its left options are ordinal numbers.
\end{Definition}

We will use small Greek letters like $\alpha$, $\beta$, $\gamma$, \dots
to denote ordinal numbers.

An ordinal number is really a number in our sense, as the following lemma
shows.

\begin{Lemma} \label{Lemma:Ordinal} \strut
      \begin{enumerate}
        \item Every ordinal number is a number.
        \item If $\alpha$ is an ordinal number, then the class of all
              ordinal numbers $\beta < \alpha$ is a set.
        \item If $\alpha$ is an ordinal number, then
              $\alpha = \game{\beta : \beta < \alpha }{\;}$, where
              $\beta$ runs through the ordinal numbers.
        \item If $\alpha$ is an ordinal number,
              then $\alpha + 1 = \game{\alpha}{\;}$.
      \end{enumerate}
\end{Lemma}
\begin{Proof}
      Note that $\beta < \alpha$ implies that there is an~$\alpha^L$ such that
      $\beta \le \alpha^L$. This follows from
      $\alpha - \beta = \game{\alpha^L - \beta}{\dots}$ and the definition
      of $G \gf 0$.
      \begin{enumerate}
        \item
          Proof by induction. By hypothesis, all $\alpha^L$ are numbers.
          Since there are no~$\alpha^R$, the condition $\alpha^L < \alpha^R$ for
          all pairs $(\alpha^L, \alpha^R)$ is trivially satisfied; hence
          $\alpha$ is itself a number.
        \item
          We claim that
          $\{\beta : \beta < \alpha\}
              = \{\alpha^L\}
                 \cup \bigcup_{\alpha^L} \{\beta : \beta < \alpha^L\}$.
          The statement follows from this by induction and the fact that
          the union of a family of sets indexed by a set is again a set.

          The RHS is certainly contained in the LHS, since all
          $\alpha^L < \alpha$.
          Now let $\beta < \alpha$. Then $\beta \le \alpha^L$ for
          some~$\alpha^L$, hence $\beta < \alpha^L$ or $\beta = \alpha^L$,
          showing that $\beta$ is an element of the~RHS.
        \item
          Let $\gamma = \game{\beta : \beta < \alpha}{\;}$ (this is a game
          by what we have just shown). Then
$\alpha-\gamma=\game{\alpha^L-\gamma}{\alpha-\gamma^L}$ (there are no
right options of $\alpha$ or $\gamma$). Since all
          $\alpha^L$ are ordinal numbers and $\alpha^L < \alpha$,
          every $\alpha^L$ is a left option of $\gamma$, hence
          $\alpha^L<\gamma$ and all $(\alpha-\gamma)^L<0$.
          By definition of $\gamma$, all $\gamma^L<\alpha$, so all
          $(\alpha-\gamma)^R>0$. Therefore $\alpha-\gamma=0$.
        \item
          We have $\alpha + 1 = \game{\alpha, \alpha^L+1}{\;}$, so we have to
          show that $\alpha^L + 1 \le \alpha$ for all~$\alpha^L$.
          Let $\beta = \alpha^L$. By induction,
          $\beta + 1 = \game{\beta}{\;} \le \game{\beta, \dots}{\;} = \alpha$.
      \end{enumerate}
\end{Proof}

Simple examples of ordinal numbers are the natural numbers:
$0 = \game{\;}{\;}$, $1 = \game{0}{\;}$,
$2 = \game{1}{\;} = \game{0, 1}{\;}$, \dots.
The next ordinal number after all the natural numbers is quite important; it is
$\omega = \game{0, 1, 2, \dots}{\;}$, the smallest infinite ordinal.

Recall that an ordered set or class is called \emph{well-ordered}
if every nonempty subset or subclass has a smallest element.
This is equivalent to the requirement that there be no infinite
descending chain of elements $x_0 > x_1 > x_2 > \dots$.

\begin{Proposition} \label{Prop:OrdWellOrdered}
      The class of ordinal numbers is well-ordered.
\end{Proposition}
\begin{Proof}
      Let $\mathcal{C}$ be some nonempty class of ordinal numbers.
      Then there is some $\alpha \in \mathcal{C}$. Replace
      $\mathcal{C}$ by the {\em set}
      $\mathcal{S} = \{\beta \in \mathcal{C} : \beta \le \alpha\}$;
      then it suffices to show that the set $\{\beta : \beta \le \alpha\}$
      is well-ordered. But every descending chain in this set is a chain
      of options of $\game{\alpha}{\;}$ and therefore must be finite by
      the~DGC.
\end{Proof}

The principle of Conway Induction applied to ordinal numbers results
in the Theorem of Ordinal Induction (sometimes called `transfinite
induction').

\begin{Theorem}[Ordinal Induction]
      Let $P$ be a property which ordinal numbers might have.
      If `$\beta$ satisfies~$P$ for all $\beta < \alpha$'
      implies `$\alpha$ satisfies $P$', then all ordinal numbers satisfy~$P$.
\end{Theorem}
\begin{Proof}
      Apply Conway Induction to the property `if $G$ is an ordinal number, then
      $G$ satisfies~$P$', and recall that $\alpha = \game{\beta < \alpha}{\;}$.
\end{Proof}

On the other hand, one could use the concept of birthdays (see below)
to prove the principle of Conway Induction from
the Theorem of Ordinal Induction.

Of course, we then also have a principle of Ordinal Recursion.
For example, we can recursively define the following numbers.

\begin{Definition}
     $2^{-\alpha} := \game{ 0 }{ 2^{-\beta} : \beta < \alpha }$
     (where $\alpha$ and $\beta$ are ordinal numbers).
\end{Definition}
These numbers are all positive and approach zero, in a similar
way as the ordinal numbers approach infinity --- for every positive
number $z > 0$, there is some ordinal number~$\alpha$ such that
$2^{-\alpha} < z$. As an example,
we have $2^{-\omega} = \omega^{-1}$. The notation is justified,
since one shows easily that $2 \cdot 2^{-(\alpha+1)} = 2^{-\alpha}$.

Finally, there is another important property of the ordinal
numbers, which is sort of dual to the well-ordering property.

\begin{Proposition} \label{Prop:OrdHasSups}
      Every \textbf{set} of ordinal numbers has a least upper
      bound within the ordinal numbers.
\end{Proposition}
\begin{Proof}
      Let $\mathcal{S}$ be such a set.
      Then $\alpha = \game{\mathcal{S}}{\;}$
      is an ordinal number\footnote{
This is the customary abuse of notation:
we mean the ordered pair $\alpha=(\mathcal{S},\{\;\})
=\game{s: s\in\mathcal{S}}{\;}$.}
and an upper bound for~$\mathcal{S}$.
      Hence the class of ordinal upper bounds is nonempty, therefore
      (because of well-ordering) there is a least upper bound.
\end{Proof}


\subsection{Birthdays}\label{Sect:Birthdays}

The concept of birthday of a game is a way of making the history
of creation of numbers and games precise. It assigns to every game
an ordinal number which can be understood as the `number of steps'
that are necessary to create this game `out of nothing' (i.e., starting
with the empty set).

\begin{Definition}[Birthday]
      Let $G$ be a game. The \emph{birthday} of~$G$, $b(G)$ is defined
      recursively by $b(G) = \game{b(G^L), b(G^R)}{\;}$.
\end{Definition}

For example, $b(0) = 0$, $b(1) = b(-1) = b(\star) = 1$; more generally, for
ordinal numbers~$\alpha$, one has $b(\alpha) = \alpha$. A game is short
if and only if its birthday is finite (see below). All non-short real
numbers have birthday~$\omega$. The successive `creation' of numbers with
the first few birthdays is illustrated by the `Australian Number Tree' in
\cite[\S~2, Fig.~2]{WW} and \cite[Fig.~0]{ONAG}.

By definition, the birthday is an ordinal number. It has the following
simple properties.

\begin{Lemma} \label{Lemma:BirthdayProperties}
      Let $G$ be a game. Then $b(G^L) < b(G)$ for all~$G^L$ and $b(G^R) < b(G)$
      for all~$G^R$. Furthermore, $b(-G) = b(G)$.
\end{Lemma}
\begin{Proof}
      Immediate from the definition.
\end{Proof}

Note that two games that are equal can have different birthdays.
For example, $\game{-1}{1} = 0$, but the first has birthday~$2$, whereas
$b(0) = 0$. But there is a well-defined minimal birthday among the
games in an equivalence class.

\begin{Proposition}
      A game~$G$ is short if and only if it has birthday $b(G) < \omega$
      (i.e., $b(G) = n$ for some $n \in \N = \{0, 1, \dots\}$).
\end{Proposition}
\begin{Proof}
      If $G$ is short, then by induction all its finitely many options
      have finite birthdays. Let $b$ be the maximum of these. Then
      $b(G) = \game{b}{\;} = b+1 < \omega$. The reverse implication follows
      from the fact that there are only finitely many games with any
      given finite birthday (this is easily seen by ordinary induction).
\end{Proof}

The birthday is sometimes useful if one needs a bound on a game.

\begin{Proposition} \label{Prop:BirthdayBound}
      If $G$ is a game, then $-b(G) \le G \le b(G)$.
\end{Proposition}
\begin{Proof}
      It suffices to prove the upper bound; the lower bound follows by replacing
      $G$ with $-G$, since both have the same birthday.

      Let $\alpha = b(G)$. $G \le \alpha$ means that for all $G^L$, we have
      $G^L \lf \alpha$, and for all $\alpha^R$, we have $G \lf \alpha^R$.
      But there are no~$\alpha^R$, so we can forget about the second condition.
      Now by induction, $G^L \le b(G^L) < b(G) = \alpha$, giving the first part.
\end{Proof}


\section{Games and Numbers}
\label{Sect:GamesNumbers}

In some sense, numbers are the simplest games --- since they are
totally ordered, we know exactly what happens (i.e., who wins) when
we add games that are numbers. It is much more difficult to deal with
general games. In order to make life easier, we try to get as much
mileage as we can out of comparing games with numbers.
References for this chapter are \cite[\S\S~8, 9]{ONAG}, \cite[\S\S~2,
6]{WW}.


\subsection{When is a game already a number?}
\label{Ssec:WhenNumber}

If we want to compare games with numbers, the first question we have
to answer is whether a given game is already (equal to) a number.
The following result gives a general recipe for deciding that two
games are equal; it can be used to provide a criterion for when
a game is a number.

\begin{Proposition}[General Simplicity Theorem]
\label{Prop:Simplicity}
      Let $G$ and~$H$ be games such that
      \begin{enumerate}
        \item $\forall G^L : G^L \lf H$ and $\forall G^R : H \lf G^R$;
        \item $\forall H^L \exists G^L : H^L \le G^L$ and
              $\forall H^R \exists G^R : H^R \ge G^R$.
      \end{enumerate}
      Then $G = H$.
\end{Proposition}
\begin{Proof}
      By the first assumption, all $G^L$ are left gift horses for~$H$, and all
      $G^R$ are right gift horses for~$H$. So $H = K$, where $K$ is the game
      whose set of left (resp.\ right) options is the union of the left
      (resp.\ right)
      options of $G$ and of~$H$. Then by the second assumption, all the options
      in~$K$ that came from~$H$ are dominated by options that came from~$G$,
      so we can eliminate all the options coming from~$H$ and get $H = K = G$.
\end{Proof}

Note that if $H$ is a number, then the second condition means that the
first condition does not hold with $H$ replaced by an option of~$H$.
This gives the usual statement of the Simplicity Theorem for comparing
games with numbers (the last claim in the following corollary).

\begin{Corollary} \label{Cor:GameIsNumber}
      Suppose $G$ is a game and $x$ is a number such that
      $\forall G^L : G^L \lf x$ and $\forall G^R : x \lf G^R$.
      Then $G$ is equal to a number; in
      fact, $G$ equals a position of~$x$.

      If no option of $x$ satisfies the assumption in place of~$x$, then $G=x$.
\end{Corollary}
\begin{Proof}
      (See also \cite[Thm.~11]{ONAG}.)
      The second statement is a special case of Prop.~\ref{Prop:Simplicity}.
      If there is an option $x'$ of~$x$ such that $G^L \lf x' \lf G^R$,
      then replace $x$ with~$x'$ and use induction.
\end{Proof}

The notion of `simple' games is defined in terms of birthdays: the earlier
a game is created (i.e.\ the smaller its birthday), the simpler it is. The
simplest game is $0=\game{\;}{\;}$ which is created first, and
subsequently more and more complicated games with later birthdays are
created out of simpler (older) ones.
The name `Simplicity Theorem' for the last statement in the corollary above
comes from the fact that in this case, $x$ is the `simplest' number
satisfying the assumption (because none of its options do). If $G$ is a
number, the  statement can then be interpreted as saying that $G$ equals
the simplest number that fits between $G$'s left and right options.

For example, a game that has no right options must equal a number,
since $G^L \le b(G^L) < b(G)$. In fact, this number is
a position of~$b(G)$, hence $G$ is even an ordinal number.

Note that a game~$G$ such that for all pairs $(G^L, G^R)$ we
have $G^L < G^R$ is not necessarily a number. A simple
counterexample is given by $G = \game{0}{\up}$, where
$\up = \game{0}{\star}$, which is a positive
game smaller than all positive (Conway) numbers (an all small game, see below
in Section~\ref{Sect:Infini}), but we have $0 < \up$.


\subsection{How to play with numbers}
\label{Ssec:PlayNumbers}

It is clear how one has to play in a number: choose an option which is as
large (or as small) as possible. It is always a disadvantage to move in a
number because one has to move to a position worse than before. In any
case, we can easily predict from the sign of the number who will win the
game, and in order to achieve this win, it is only necessary to choose
options of the correct sign. As far as playing is concerned, numbers are
pretty boring! But what is good play in a sum $G + x$, where $x$ is a
number and $G$ is not?

\begin{Theorem}[Weak Number Avoidance Theorem]
      \label{Thm:WeakNumberAvoidance}
      If $G$ is a game that is not equal to a number and $x$ is a number, then
      \[ x \lf G \iff \exists G^L : x \le G^L \,. \]
\end{Theorem}
\begin{Proof}
      The implication `$\Leftarrow$' is trivial. Now assume $x \lf G$.
      This means that either $\exists G^L \ge x$ (and we are done),
      or $\exists x^R \le G$, and we can assume that for all~$G^L$, we have
      $x \gf G^L$. Since $G$ is not equal to a number by assumption,
      Cor.~\ref{Cor:GameIsNumber} implies that there is some~$G^R$ with
      $x \ge G^R$. But then we get $G \ge x^R > x \ge G^R$, in contradiction
      to the basic fact $G \lf G^R$.
\end{Proof}

If we apply this to $G$ and~$-x$, it says
$G + x \gf 0 \iff G^L + x \ge 0$ for some~$G^L$. In words, this
means that if there is a winning move in the sum $G + x$, then
there is already a winning move in the $G$~component. In short:
\begin{quote}
      \em In order to win a game, you do not have to move in a
      number, unless there is nothing else to do.
\end{quote}
This does {\em not}
mean, however, that the other options $G + x^L$ are redundant (i.e.,
dominated or reversible).
This is only the case in general when $G$ is {\em short}.
In order to prove this stronger version of the Number Avoidance
Theorem, we need some preparations.

\begin{Definition}[Left and Right Stops] \label{DefLeftRightStops}
      Let $G$ be a short game. We define (short) numbers
      $L(G)$ and~$R(G)$, the \emph{left} and \emph{right stops}
      of~$G$, as follows.

      If $G$ is (equal to) a number, we set $L(G) = R(G) = G$. Otherwise,
      \[ L(G) = \max_{G^L} R(G^L) \quad\text{and}\quad
         R(G) = \min_{G^R} L(G^R) \,.
      \]
\end{Definition}

Since $G$ has only finitely many options, the maxima and minima
exist. Note that $G$ must have left and right options; otherwise
$G$ would be a number.

One can think of~$L(G)$ as the best value Left can achieve as
first player in~$G$, at the point when the game becomes a number.
Similarly, $R(G)$ is the best value Right can achieve when moving
first. (Since numbers are pretty uninteresting games, we can stop
playing as soon as the game turns into a number; this number can then
conveniently be interpreted as the score of the game.)

\begin{Proposition} \label{Prop:LRStops}
      Let $G$ be a \textbf{short} game.
      \begin{enumerate}
        \item If $y$ is a number, then the following implications hold.
              \begin{align*}
                 y > L(G) \Longrightarrow y > G\,, &\qquad
                 y < L(G) \Longrightarrow y \lf G\,, \\
                 y < R(G) \Longrightarrow y < G\,, &\qquad
                 y > R(G) \Longrightarrow y \gf G\,.
              \end{align*}
        \item If $z > 0$ is a positive number and $G$ is not equal to a number,
              then $G < G^L + z$ for some~$G^L$.
      \end{enumerate}
\end{Proposition}
\begin{Proof}
      \begin{enumerate}
        \item
          If $G$ is a number, then $G = L(G) = R(G)$, and the statements
          are trivially true. If $G$ is not a number, we proceed by
          induction.

          Assume $y > L(G)$. Then by definition, we have
          $y > R(G^L)$ for all~$G^L$. By induction hypothesis, this
          implies $y \gf G^L$ for all~$G^L$.
          By Thm.~\ref{Thm:WeakNumberAvoidance}, this means $y \ge G$,
          hence $y > G$, since $y$ cannot equal~$G$ ($y$ is a number,
          but $G$ is not).

          Now assume $y < L(G)$. Then by definition, $y < R(G^L)$ for
          some~$G^L$. By induction hypothesis, $y < G^L$ for this~$G^L$.
          But this implies $y \lf G$.

          The other two statements are proved analogously.
        \item
          Let $y = L(G) + z/2$. By the first part, we then have $y > G$ and
          $y-z \lf G$. By Thm.~\ref{Thm:WeakNumberAvoidance}, there is
          a~$G^L$ with $y-z \le G^L$. Hence $G < y \le G^L + z$.
      \end{enumerate}
\end{Proof}

The first part of the preceding proposition can be interpreted as
saying that the {\em confusion interval} of~$G$ (the set of numbers
$G$ is fuzzy to) extends from $R(G)$ to~$L(G)$. (The endpoints may
or may not be included; this depends on who has the move when the
game reaches its stopping position; compare the Temperature Theory
of games in~\cite[\S~9]{ONAG} and \cite[\S~6]{WW}.)

With these preparations, we can now state and prove the Number
Avoidance Theorem in its strong form.

\begin{Theorem}[Strong Number Avoidance Theorem]
      If $G$ is a \textbf{short} game that is not equal to a number and
      $x$ is a number, then $G + x = \game{ G^L + x }{ G^R + x }$.
\end{Theorem}
\begin{Proof}
      (Compare \cite[Thm.~90]{ONAG}.)\\
      We have $G + x = \game{ G^L + x, G + x^L }{ G^R + x, G + x^R }$. Consider
      an option $G + x^L$. By the second part of Prop.~\ref{Prop:LRStops},
      applied to $z = x-x^L$,
      there is some~$G^L$ such that $G < G^L + x - x^L$. Therefore
      $G + x^L < G^L + x$ is dominated and can be removed. An analogous
      argument applies to $G + x^R$.
\end{Proof}

Let us demonstrate that the assumption on~$G$ is really necessary.
Consider the game
\[ G = \game{\Z}{\Z}
        = \game{\dots,-2,-1,0,1,2,\dots}{\dots,-2,-1,0,1,2,\dots} \,\,.
\]
Then $\game{G^L + 1}{G^R + 1} = \game{\Z + 1}{\Z + 1} = \game{\Z}{\Z} = G$,
but (of course) $G + 1 \neq G$, since $1 \neq 0$.

This game~$G$ also provides a counterexample to the second part
of Prop.~\ref{Prop:LRStops} for non-short games. It is easily seen that $G
\fu n$ for all integers~$n$, hence $G < G^L + 1$ is impossible.

The deeper reason for this failure is that it is not possible to
define left and right stops for general games. And this is because
there is nothing like a supremum of an arbitrary set of numbers
(for example, $\Z$ has no least upper bound). This should be
contrasted with the situation we have with ordinal numbers, where
every set of ordinal numbers has a least upper bound within the
ordinal numbers.


\section{Infinitesimal games \cite[\S~8]{WW}} \label{Sect:Infini}


If a game is approximately the size of a positive (or negative) number,
we know that Left (or Right) will win it. But there are games which are
less than all positive numbers and greater than all negative numbers,
and we do not get a hint as to who is favored by the game. Such a game is
called {\em infinitesimal}.

\begin{Definition}[Infinitesimal] \strut
      \begin{enumerate}
        \item A game $G$ is called \emph{infinitesimal} if
              $-2^{-n} < G < 2^{-n}$ for all natural numbers~$n$.
        \item A game $G$ is called \emph{strongly infinitesimal}
              if $-z < G < z$ for all positive numbers~$z$.
      \end{enumerate}
\end{Definition}

An example of an infinitesimal, but not strongly infinitesimal
game is given by $2^{-\omega} = \game{ 0 }{ 2^{-n} : n \in \N }$.
The standard example of a positive strongly infinitesimal game
is $\up = \game{0}{\star}$ (pronounced `up').
More examples are provided by the following class of games.

\begin{Definition}[All Small]
      A game~$G$ is called \emph{all small} if every position of~$G$
      that has left options also has right options and vice versa.
\end{Definition}

Since a game without left (or right) options is equal to a number (see the example
after Corollary~\ref{Cor:GameIsNumber}), this definition is equivalent to `no number
other than~$0$ occurs as a position of~$G$'. The simplest all small games are $0$,
$\star$, $\up = \game{0}{\star}$ and $\down = \game{\star}{0}$ (note that
$\game{\star}{\star} = 0$). They show that an all small game can fall into
any of the four outcome classes.

\begin{Proposition} \label{Prop:AllSmall}
      If $G$ is an all small game, then $G$ is strongly infinitesimal.
\end{Proposition}
\begin{Proof}
      If $G$ is a number, then $G = 0$, and the claim holds trivially.
      Otherwise, $G$ has left and right options, all of which are
      strongly infinitesimal by induction. Let $z > 0$ be a number.
      Let $z^R$ be some right option of~$z$.
      Then for all and hence for some~$G^R$, we have $G^R < z < z^R$ and
      so $G \lf z^R$. On the other hand, we have $G^L < z$ for all~$G^L$.
      Together, these two facts imply that $G \le z$. The inequality
      $G \ge -z$ is shown in the same way.
\end{Proof}

Not all strongly infinitesimal games are all small. Examples
are provided by `tinies' and `minies' like $\game{0}{\game{0}{-1}}$
\cite[\S~5]{WW}.

We will see in a moment that for short games,
`infinitesimal' and `strongly infinitesimal' are the same.
More precisely, the following theorem tells us that a short infinitesimal
game is already bounded by some integral multiple of~$\up$.
(This is mentioned in \cite[\S~20: `The Paradox']{WW}, but we do not
know of a published proof.)
For example, we have $\star < 2 \up$, and so also $2 \up + \star > 0$.

\begin{Theorem} \label{Thm:ShortInfinitesimal}
      Let $G$ be a \textbf{short} game.
      \begin{enumerate}
        \item If $G \lf 2^{-n}$ for all positive integers~$n$, there
              is some positive integer~$m$ such that $G \lf m \up$.
        \item If $G \le 2^{-n}$ for all positive integers~$n$, there
              is some positive integer~$m$ such that $G \le m \up$.
      \end{enumerate}
\end{Theorem}
\begin{Proof}
      We can assume that $G$ is not a number, because otherwise the
      assumption implies in both cases that $G \le 0$. We use induction.
      \begin{enumerate}
        \item
          $G \lf 2^{-n}$ means by Thm.~\ref{Thm:WeakNumberAvoidance}
          that there is some~$G^R$ with $G^R \le 2^{-n}$. Since there
          are only finitely many~$G^R$ ($G$ is short), there must be
          one~$G^R$ that works for infinitely many and hence for all~$n$.
          By induction, we conclude that $G^R \le m \up$ for this~$G^R$
          and some~$m$ and therefore $G \lf m \up$ for this~$m$.
        \item
          \begin{itemize}
            \item[(a)]
              If $G \le 2^{-n}$, then $G < 2^{-(n-1)}$. Hence $G$ satisfies
              the assumption of the first part, and we conclude that there
              is some $m_0$ such that $G \lf m \up$ for all $m \ge m_0$.
            \item[(b)]
              $G \le 2^{-n}$ implies that for all~$G^L$, we have $G^L
\lf 2^{-n}$.
              By induction, for every~$G^L$, there is some $m$ such that
              $G^L \lf m \up$. Since there are only finitely many~$G^L$
              ($G$ is short), there is some $m_1$ such that $G^L \lf m_1 \up$
              for all~$G^L$.
            \item[(c)]
              Now let $m = \max\{m_1, m_0+3\}$. By~(a), we know that
              $G \lf (m-3) \up < (m-1) \up + \star = (m\up)^R$. By~(b), we
              know that $G^L \lf m \up$ for all~$G^L$. Together, these
              imply $G \le m \up$.
          \end{itemize}
      \end{enumerate}
\end{Proof}

The apparent asymmetry in this proof is due to the lack of something
like an `Up Avoidance Theorem'.

This result says that infinitesimal short games can be measured in
(short) units of~$\up$. This is the justification behind the `atomic weight
calculus' described in~\cite[\S\S~7, 8]{WW}.

\begin{Remark}
      It is perhaps tempting to think that a similar statement should
      be true for strongly infinitesimal general games. If one tries
      to mimic the above proof, one runs into two difficulties.
      In the first part, we have used that any finite set of positive
      short numbers has a positive short lower bound. The corresponding
      conclusion would still be valid, since any set of positive numbers
      has a positive lower bound (which is a number). (For $z > 0$ it
      is easy to see that $z \ge 2^{-b(z)}$; the claim then follows
      from the statement on upper bounds for sets of ordinals.)
      In the second part, we have used that any finite set of natural
      numbers (or of multiples $m \up$) has an upper bound. This does not
      seem to generalize easily.
\end{Remark}


\section{Impartial Games \cite[\S~3]{WW}, \cite[\S~11]{ONAG}}
\label{Sect:Impartial}


\subsection{What is an impartial game?}
\label{Ssect:Impartial}

An impartial game is one in which both players have the same possible
moves in every position. Formally, this reads as follows.

\begin{Definition}[Impartial Game]
\label{DefImpartial}
An {\em impartial game} is a game for which the sets of left and right
options are equal, and all options are impartial games themselves.
\end{Definition}

It follows that every impartial game $G$ satisfies $G=-G$, hence $G+G=0$.
Therefore, $G=0$ or $G\fu 0$. There is no need to distinguish the sets of
left and right options, so we simply write $G=\{G',G'',\dots\}$, where
$\{G',G'',\dots\}$ is the set of options of $G$ (again, this notation is
not meant to suggest that the set of options should be countable or
non-empty).

The standard examples include the game of Nim: it consists of a finite
collection of heaps $H_1,\dots,H_n$, each of which is an ordinal number
(some number of coins, matches, etc.; if this makes you feel more
comfortable, think
of natural numbers only). A move consists of reducing any
single heap by an arbitrary amount (i.e., replacing one of the ordinal
numbers by a strictly smaller ordinal number), leaving all other heaps
unaffected. A move in the ordinal number $0$ is of course impossible, so
the game ends when all heaps are reduced to $0$. As usual, winner is the
one who made the last move. Note that a single heap game is trivial: if
the heap is non-zero, then the winning move consists in reducing the heap
to zero, leaving no legal move to the opponent. A game with several heaps
is the sum (in our usual sense) of its heaps: it is $H_1+H_2+\dots+H_n$.

Conway coined the term {\em nimber} for a single Nim heap, and he writes
$\star n$ for a heap of size $n$ (where $n$ is of course an ordinal number).
The rules of Nim can then simply be written as $\star 0 = 0$,
$\star n = \{\star 0, \star 1, \dots, \star(n-1)\}$ (if $n$ is finite) or
$\star n = \{\star k : k < n\}$
(in general).
We have $\star n = \star k$ if and only $n=k$:
the equality $\star n = \star k$ means
$0 = \star n - \star k = \star n + \star k$ (note that $\star k = -\star k$
since nimbers are impartial games), and in $\star n + \star k$ with $n \neq k$
the first player wins by reducing the larger
heap so as to leave two heaps of equal size to the opponent). The nimbers
inherit a total ordering from the ordinal numbers so that every set of
nimbers is well-ordered (Proposition~\ref{Prop:OrdWellOrdered}).
But note that this is not the same as the ordering of general games
restricted to impartial games: a Nim heap $\star n$ of size $n>0$ has
$\star n \fu 0$, not $\star n > 0$!

\medskip
{\bf A second note on set theory.}
Since there is no need to distinguish the left and right options of impartial
games, and all these options are impartial games themselves,
Definition~\ref{DefGame} of Games simplifies to the following:

\begin{Definition}[Impartial Game] \label{DefImpartialGame}
\strut
     \begin{enumerate}
       \item \label{DefGameImp1}
         Let $G$ be a set of impartial games. Then
         $G$ is an {\em impartial game.}
       \item \label{DefGameImp2}
         (Descending Game Condition for Impartial Games). There is no infinite
sequence of games
         $G^i$ with $G^{i+1}\in G^i$ for all $i\in\N$.
     \end{enumerate}
\end{Definition}
\nopagebreak

Therefore, impartial games are just sets in
the sense of Zermelo and Fraenkel; the Descending Game Condition exactly
reduces to the Axiom of Foundation.\footnote{The
observation that every set can be viewed as an impartial game leads to
amusing questions of the type ``what is the Nim heap equivalent to $\Z$,
$\Q$, $\R$, $\C$, \dots?'' The answer of course depends on the exact way of
representing these sets within set theory.}



\subsection{Classification of impartial games}
\label{Ssect:SpragueGrundy}

There is a well-known classification of impartial games, due to Sprague
and Grundy, which says that every impartial game is equal to a
well-defined nimber $\star n$, hence equal to a Nim game with a single heap,
which has trivial winning strategy. Unfortunately, this does not make
every impartial game easy to analyze: in practice, it might be hard to
tell exactly which nimber an impartial game is equal to; as an example, we
mention the game of {\em Sylver Coinage}\footnote{Sylver Coinage is
usually played in mis\`ere play; however, one can equivalently declare
the number $1$ illegal and use the normal winning convention.}
\cite[\S~18]{WW}.

The theory of impartial games is based on the following definition.
\begin{Definition}[The mex: \emph{m}inimal \emph{ex}cluded nimber]
\label{DefMex}
     Let $G$ be a set of nimbers. Then $\mex(G)$ is the least nimber not
     contained in $G$.
\end{Definition}

Note that every {\em set} of nimbers has an upper bound by
Proposition~\ref{Prop:OrdHasSups}, and
the well-ordering of nimbers implies that every {\em set} of nimbers has a
well-defined minimal excluded nimber, so the $\mex$ is well-defined.

The first and fundamental step of the classification is the following
observation.

\begin{Proposition}[Bogus Nim]
\label{PropBogusNim}
Let $G$ be any set of nimbers. Then $G=\mex(G)$.
\end{Proposition}
Some remarks might be in order here. First, $G$ is of course an impartial
game itself, describing its set of options which are all nimbers. Now
$\mex(G)$ is a particular nimber, and the result asserts that the game $G$
is equal as a game to the nimber $\mex(G)$.

\begin{Proof}
All we need to show is that $G-\mex(G)=G+\mex(G)=0$. The main trick is to write
this right: $\mex(G)=\star n=\{\star k : k<n\}$ for an ordinal number $n$,
while $G=\star n\cup G'$ where $G'$ is the set of options of $G$ exceeding
$\star n$ (note that by assumption $\star n$ is not an option of $G$). The
first player has three kinds of possible moves: move in $\mex(G)=\star n$ to
some $\star k<\star n$, move in $G$ to an option $\star k<\star n$, or move in
$G$ to an option $\star k>\star n$. The first leads to $G+\star k$ and is
countered by the move in $G$ to $\star k+\star k=0$; the second kind leads to
$\star k+\star n$ and is countered by the move in $\star n$ to $\star k+\star
k=0$; finally, the third leads to $\star k+\star n$ (this time, with $k>n$) and
is countered by a move in $\star k$ to $\star n+\star n=0$. In all three cases,
the second player moves to $0$ and wins.
\end{Proof}

The name `Bogus Nim' refers to the following interpretation of this game:
the game $G$ really is a Nim heap $\star n$ (offering moves to all
$\star k$ with $k < n$), but in addition it is allowed to increase the
size  of the heap.
This increasing is immediately reversed by the second player, bringing the
heap back to $\star n$ (in which no further increase is possible), so all
increasing moves are reversible moves.

A rather obvious corollary is the following: $G=\star 0$ if and only if no
legal move in $G$ leads to $\star 0$ (if $\star 0$ is no option of
$G$, then clearly
$\mex(G)=\star 0$): if $G$ has the option $\star 0$, then this is a
winning move for
the first player, hence $G\fu 0$; otherwise, all options of $G$ (if any)
lead to nimbers $\star n\neq \star 0$ from which the second player wins.

\begin{Theorem}[The Classification of Impartial Games]
\label{ThmSpragueGrundy} \lineclear
Every impartial game is equal to a unique nimber.
\end{Theorem}
\begin{Proof}
We use Conway induction: write $G=\{G',G'',\dots\}$, listing all options
of $G$. By the inductive hypothesis, every option of $G$ is equal to a
nimber, hence (using Theorem~\ref{ThmEqualGames})
\[
G=\{\star k : \text{ there is an option of $G$ which equals $\star k$}\}
\,\,.
\]
Therefore, all options of $G$ are equal to nimbers, hence by
Proposition~\ref{PropBogusNim} $G$ is equal to the $\mex$ of all its
options, which is a nimber. Uniqueness is clear.
\end{Proof}

We should mention that for the game Nim, there is a well-known explicit
strategy, at least for heaps of finite size: in the game
$G=\star n_1+\star n_2+\dots+\star n_s$, write each heap size $n_i$
in binary form and
form the {\sc exclusive or (XOR)} of them. Then $G=0$ if and only if the
{\sc XOR} is zero in every bit: in the latter case, it is easy to see that
every option of $G$ is non-zero, while if the {\sc XOR} is non-zero, then
every heap which contributes a {\tt 1} to the most significant bit of the
{\sc XOR} can be reduced in size so as to turn the game into a zero game.

Nim is often played in Mis\`ere play, and the binary strategy as just
described works in this form almost without difference, except very near
the end when there are at most two heaps of size exceeding one. It is
sometimes wrongly concluded that there was a general theory of impartial
games in Mis\`ere play similarly as Theorem~\ref{ThmSpragueGrundy}.
However, this is false: the essential difference is the use of
Theorem~\ref{ThmEqualGames} which allows us to replace an option by an
equivalent one, and this is based on the usual winning convention. There
is no analog in Mis\`ere play to the Sprague-Grundy-Theory: for any two
impartial games $G$ and $H$ without reversible options which are different
in form (i.e., $G\not\equiv H$), there is another impartial game $K$ such
that the winners of $G+K$ and $H+K$ are different; see
\cite[Chapter~12]{ONAG}.




\end{document}